\documentclass{scrartcl}

\makeatletter
\DeclareOldFontCommand{\rm}{\normalfont\rmfamily}{\mathrm}
\DeclareOldFontCommand{\sf}{\normalfont\sffamily}{\mathsf}
\DeclareOldFontCommand{\tt}{\normalfont\ttfamily}{\mathtt}
\DeclareOldFontCommand{\bf}{\normalfont\bfseries}{\mathbf}
\DeclareOldFontCommand{\it}{\normalfont\itshape}{\mathit}
\DeclareOldFontCommand{\sl}{\normalfont\slshape}{\@nomath\sl}
\DeclareOldFontCommand{\sc}{\normalfont\scshape}{\@nomath\sc}
\makeatother

\usepackage{arxiv}
\usepackage{authblk}

\usepackage{listings}
\usepackage[frozencache]{minted}
\usepackage{subcaption}
\usepackage{paralist}
\newcommand{\pder}[2]{\frac{\partial#1}{\partial#2}}

\newcommand{\vertiii}[1]{{\left\vert\kern-0.25ex\left\vert\kern-0.25ex\left\vert #1 
    \right\vert\kern-0.25ex\right\vert\kern-0.25ex\right\vert}}

\usepackage{tikz}
\usepackage{pgfplots}
\pgfplotsset{compat=1.9}
\pgfkeys{/pgfplots/tuftelike/.style={
  semithick,
  tick style={major tick length=4pt,semithick,black},
  separate axis lines,
  axis x line*=bottom,
  axis x line shift=10pt,
  axis y line*=left,
  axis y line shift=10pt}}
  
\usetikzlibrary{pgfplots.dateplot}
\usepackage{pgfplotstable}
\usepackage{pgf}
\usetikzlibrary{fit, positioning, shapes, shadows, calc, intersections, through, backgrounds, arrows, matrix, patterns}

\pgfdeclarelayer{bg}    %
\pgfdeclarelayer{fg}    %
\pgfsetlayers{bg,main,fg}  %

\definecolor{participantAcolor}{RGB}{243,98,33}
\definecolor{participantBcolor}{RGB}{0,102,189}
\definecolor{participantCcolor}{RGB}{66,70,50}
\definecolor{cplColor}{RGB}{0,0,0}
\definecolor{preciceColor}{RGB}{ 161, 177, 25}

\makeatletter
\pgfkeys{%
  /tikz/on layer/.code={
    \pgfonlayer{#1}\begingroup
    \aftergroup\endpgfonlayer
    \aftergroup\endgroup
  },
  /tikz/node on layer/.code={
    \gdef\node@@on@layer{%
      \setbox\tikz@tempbox=\hbox\bgroup\pgfonlayer{#1}\unhbox\tikz@tempbox\endpgfonlayer\egroup}
    \aftergroup\node@on@layer
  },
  /tikz/end node on layer/.code={
    \endpgfonlayer\endgroup\endgroup
  }
}
\def\node@on@layer{\aftergroup\node@@on@layer}
\makeatother

\usepackage{amsmath}
\usepackage{amssymb}
\usepackage{hyperref}
\usepackage[perpage]{footmisc}
\usepackage{booktabs}

\usepackage{mathtools}
\let\oldtextbf\textbf
\renewcommand{\textbf}[1]{\oldtextbf{\boldmath #1}}

\newcommand{\DSC}{\mathcal{D}_{\text{SC}}}
\newcommand{\NSC}{\mathcal{N}_{\text{SC}}}
\newcommand{\DWI}{\mathcal{D}_{\text{WI}}}
\newcommand{\NWI}{\mathcal{N}_{\text{WI}}}
\newcommand{\DS}{\mathcal{D}}
\newcommand{\NS}{\mathcal{N}}
\newcommand{\nD}{n_\mathcal{D}}
\newcommand{\nN}{n_\mathcal{N}}
\newcommand{\pD}{p_\mathcal{D}}
\newcommand{\pN}{p_\mathcal{N}}
\newcommand{\cD}{{c_\mathcal{D}}}
\newcommand{\cN}{{c_\mathcal{N}}}
\newcommand{\cDe}{\cD_{\text{end}}}
\newcommand{\cNe}{\cN_{\text{end}}}
\newcommand{\I}{\mathcal{I}}

\newcommand{\gpolorder}{\alpha}

\usepackage{cleveref}

\title{Quasi-Newton Waveform Iteration for Partitioned Fluid-Structure Interaction}

\author[1]{Benjamin R\"uth}
\author[2]{Benjamin Uekermann}
\author[3]{Miriam Mehl}
\author[4]{Philipp Birken}
\author[4]{Azahar Monge}
\author[1]{Hans-Joachim Bungartz}

\affil[1]{Department of Informatics, Technical University of Munich, Germany, \texttt{\{rueth,bungartz\}@in.tum.de}}
\affil[2]{Department of Mechanical Engineering, Einhoven University of Technology, The Netherlands, \texttt{b.w.uekermann@tue.nl}}
\affil[3]{Institute for Parallel and Distributed Systems, University of Stuttgart, Germany, \texttt{miriam.mehl@ipvs.uni-stuttgart.de}}
\affil[4]{Centre for the Mathematical Sciences, Lund University, Sweden, \texttt{philipp.birken@math.lu.se}, \texttt{azahar.sz@gmail.com}}

\begin{document} 

\maketitle

\begin{abstract}
We present novel coupling schemes for partitioned multi-physics simulation that combine four important aspects for strongly coupled problems: implicit coupling per time step, fast and robust acceleration of the corresponding iterative coupling, support for multi-rate time stepping, and higher-order convergence in time. 
To achieve this, we combine waveform relaxation -- a known method to achieve higher order in applications with split time stepping based on continuous representations of coupling variables in time -- 
with interface quasi-Newton coupling, which has been developed throughout the last decade and is generally accepted as a very robust iterative coupling method even for gluing together black-box simulation codes. 
We show convergence results (in terms of convergence of the iterative solver and in terms of approximation order in time) for two academic test cases -- a heat transfer scenario 
and a fluid-structure interaction simulation. We show that we achieve the expected approximation order and that our iterative method is competitive in terms of iteration counts with those designed for simpler first-order-in-time coupling.
\end{abstract}


\section{Introduction}\label{sec::intro}

We consider time dependent fluid-structure interaction (FSI) problems. Examples are numerous: flutter of airplanes or wind turbines, blood flow in arteries, gas quenching or cooling of rocket engines, to only name a few. We follow a partitioned coupling approach, which abides the natural wish to re-use established single-physics solvers. Hereby, our goal is to have highest flexibility, while still having a fast time to solution in a black-box coupling framework. Only minimal information from the solvers is used in that only boundary conditions at the interface are exchanged, as well as geometric interface information, that can be used to map non-matching interface meshes. 

In this article, we focus on the time integration. The combination of different physics oftentimes implies different time scales. Thus, we would like to be able to use different time steps in the involved solvers. Additionally, higher-order time integration is required to efficiently produce sufficiently accurate results. An open problem is to design stable, efficient and high-order methods that allow for different time steps in the separate solvers \cite{Keyes2013} and for the use of black-box solvers. A class of methods that potentially give these properties are waveform relaxations, also called waveform or dynamic iterations \cite{Gander2015, Burrage1995}. These integrate over a time window, allowing the subsolvers to use potentially arbitrary methods and time steps in that window, in the spirit of the partitioned approach. To this end, continuous interpolations in time are employed to provide solution values to the other problem. High order can be obtained by sufficiently high order in the subsolvers, as well as the interpolation. To get stability, the time integration over the window is repeated iteratively. An efficient method is obtained, when this iteration converges fast. 

The standard option is to add a relaxation step, which is why the terms waveform relaxation and waveform iteration are used synonymously in most literature.
Waveform relaxation is widely applied in field/circuit coupling \cite{Schops2011, Schops2012, Maciejewski2017}. 
Recently, Neumann-Neumann waveform relaxation was suggested \cite{gander2016dirichlet} and applied to conjugate heat transfer problems \cite{Monge2018, Monge2018Thesis}. Optimal relaxation parameters are determined for this specific problem, leading to superlinear convergence. This was later extended to the time adaptive case \cite{Monge2019_Adaptive}.
The optimal relaxation parameter is, however, problem dependent. This means that when applied as a black-box method, many iterations are needed.

Different alternatives have been suggested to speed up waveform relaxation, e.g. convolution waveform methods \cite{rewhal:95}, Krylov methods \cite{lumwu:03} and multigrid \cite{horvan:95}. Another alternative are Waveform relaxation Newton methods \cite{whodvr:85}. There, the Newton method is formulated on the initial value problems using Fr\'{e}chet derivatives. All of these offer interesting theoretical properties, which are however, hard to exploit in practice.

Another idea that is more straightforward to use is to use a quasi-Newton method instead of relaxation to improve the convergence of waveform iteration.
This has been very successfully applied for non-waveform partitioned black-box FSI,
and, due to its efficiency, is now state of the art in this field. The basic variant is called Interface Quasi-Newton Inverse Least-Squares (IQN-ILS) \cite{Degroote2009_Performance, Haelterman2009_QNLS}. As a non-linear solver, the method is also known as Anderson acceleration \cite{anders:65, fansaa:09}. Important recent improvements have been parallel execution of solvers \cite{Mehl2016}, filtering \cite{Haelterman2015_Filtering}, reducing tuning parameters \cite{Bogaers2014_MVQN, Lindner2015_MVQN, Scheuf:QN} and the generalization to multiple coupled solvers \cite{Bungartz2014_MIQN_CompuMech}. The resulting methods lead to very few coupling iterations, are robust without tuning parameters, and have linear compute and memory complexity in terms of interface degrees of freedom due to matrix-free formulations and due to a mesh independent convergence rate. Implicit coupling together with state-of-the-art IQN methods overcome added mass instability and convergence issues and can therefore also be applied for bio-mechanical applications\cite{Davis2019, Uekermann2016}. 

In this paper, we combine waveform iteration (WI) with IQN, with the goal of getting a novel fast partitioned black-box coupling method that still has high order and allows for different time step sizes in different solvers. To do this, we consider versions of the quasi-Newton method that differ in terms of which subsets of time step data in the so-called coupling windows are used to approximate the system Jacobian.
A similar approach has been followed in \cite{DeMoerloose2019} for a less general setting. There, subcycling, i.e., using several time steps of one solver within one time step of the other one, is considered. For a 1D problem and certain time step ratios, the stability of the time integration is analyzed, which is confirmed by numerical results in 2D. Our approach allows different numbers of time steps of both solvers within a so-called coupling window such that the time steps of involved solvers no longer have to be multiples of each other (called multi-rate time stepping).

For the description and the numerical experiments, we only consider so called Dirichlet-Neumann type coupling methods. These are a basic building block in FSI that only use boundary conditions that are available in the subsolvers anyhow. To this end, Dirichlet, respectively Neumann data at the interface are exchanged between the solvers. A generalization for other coupling conditions is possible, but requires additional numerical testcases to evaluate the resulting performance. 
We do not adapt time steps dynamically in time, but we would like to point out that time adaptivity is straightforward to use in waveform iterations. This is important, since it can give dramatic performance improvements \cite{biquhm:11,biquhm:10}. 

We demonstrate up to third order in time for the combined method. The iteration numbers do not vary much when varying the time step ratios and are rather small. 
We test the proposed algortihms by using the coupling library preCICE\cite{Bungartz2016_preCICE} with the finite element frameworks FEniCS\cite{Logg2012} and Nutils\cite{Nutils}. We, therefore, use IQN from preCICE and implement WI in a new middle layer between preCICE and the coupled solvers.

The structure of the article is as follows: In Sect.~\ref{sec::coupling}, we explain the mathematical framework for both concepts -- waveform iteration for multi-rate coupled time stepping on the one hand and quasi-Newton on the other hand -- separately and then show how we can combine them. We then describe in Sect.~\ref{sec::software} how we realized a prototype implementation in a middle software layer between coupled solvers and preCICE. In Sect.~\ref{sec::results}, we show numerical results for a coupled heat equation and an FSI scenario, both in 2D. We conclude in Sect.~\ref{sec::conclusions}. 

\section{Coupling Algorithms}\label{sec::coupling}

Before introducing the implementation in the coupling library preCICE, we present the details of our coupling approach allowing for multi-rate time stepping, higher-order accuracy in time, and fast iterative convergence in several steps:
Section \ref{sec:basic} introduces our notation and the general form of the coupling formulation as used in standard low-order approaches.
Section \ref{sec:waveform} explains how to extend this basic coupling with a waveform approach, which includes information from individual time steps to increase the approximation order. Section \ref{sec::quasi-newton} explains how to construct a quasi-Newton scheme for a general setting. Section \ref{sec::qn-wc} combines both concepts to define quasi-Newton variants for multi-rate higher-order coupling.

\subsection{Basic Notation and Multi-rate Coupling via Single Value Coupling}
\label{sec:basic}

We consider a coupled problem on two non-overlapping domains $\Omega_{\mathcal{D}}$ and $\Omega_{\mathcal{N}}$.
We use a partitioned approach, i.e., we rely on stand-alone solvers for each domain, which we regard as black boxes providing access only to input and output data, but not to discretization details.
The coupling of solvers happens at the coupling interface $\Gamma = \Omega_{\mathcal{D}} \cap \Omega_{\mathcal{N}}$. 
We couple the solvers via Dirichlet-Neumann coupling. The Dirichlet solver on $\Omega_{\mathcal{D}}$ takes a Dirichlet boundary condition at the coupling interface $\Gamma$ from the Neumann solver on $\Omega_{\mathcal{N}}$, while the Neumann solver takes a Neumann boundary condition at the coupling interface $\Gamma$ from the Dirichlet solver.
For FSI, the natural choice is to use the fluid solver as Dirichlet solver with Dirchlet boundary values being velocities or displacements, and the solid solver as Neumann solver with Neumann boundary values being forces or stresses (see Sect. \ref{ssec:FSI}). For other coupled problems such as heat transfer, which we use as a second application test case in Sect.~\ref{ssec:heat}, this choice has to be made based on material parameters of the two media \cite{monbir:18}. Here, Dirichlet boundary values represent temperature values at the coupling interface and Neumann boundary values the heat flux computed from the internal temperature field.

Both solvers perform time stepping (explicit or implicit) using their respective schemes in a given time window $[t_{\text{ini}}; t_{\text{ini}} + \Delta t]$ with initial time $t_{\text{ini}}$ and window size $\Delta t$. Within each window, we assume them to use a fixed, but possibly different time step size. The Dirichlet and the Neumann solver, thus, perform $\nD$ or $\nN$ time steps of size $\delta t_\mathcal{D}$ and $\delta t_\mathcal{N}$, respectively ($\Delta t = \nD \delta t_\mathcal{D} = \nN \delta t_\mathcal{N}$).

We now define the operators $\DSC$ and $\NSC$ as the action of our black-box solvers over the complete window, given Dirichlet boundary values $\cDe$ at the end of the coupling window (and initial values at the boundary and inside the domain, which we omit as arguments for better readability), the Dirichlet solver returns Neumann boundary values $\cNe$ at the end of the coupling window. The latter are used as an input for the Neumann solver that returns Dirichlet boundary values $\cDe$, again:
\[ \DSC: \cDe \mapsto \cNe \; \mbox{ and } \; \NSC: \cNe \mapsto \cDe. \]
The coupling condition at the interface can now be expressed by means of the fixed-point equation
\begin{equation}
\cDe = \underbrace{\NSC \circ \DSC}_{\mbox{$=:H_{\text{SC}}$}} ( \cDe ) \label{eq:FPSC}
\end{equation}
acting only on the Dirichlet interface data $\cDe$ at the end of the window as the two codes exchange only bounday data at this single point in time.  $H_{\text{SC}}$ denotes the corresponding fixed-point operator. We refer to coupling methods solving Eq.~(\ref{eq:FPSC}) in this contribution as single value coupling and denote them as \textbf{SC($\nD,\nN$)}. Further, the execution of $\nD$ and $\nN$ time steps in the solvers before coupling at the end of time window is called multi-rate time stepping. the special case with either $\nD = 1, \nN \neq 1$ or $\nD \neq 1, \nN = 1$ is commonly called subcycling\cite{Farhat2000, Felippa2004}.

Eq.~(\ref{eq:FPSC}) entails a Gauss-Seidel-type coupling, executing $\DSC$ and $\NSC$ one after the other, while the roles of $\DSC$ and $\NSC$ could also be swapped leading then to a fixed-point equation in terms of $\cN$. Furthermore, we could instead also consider a Jacobi-type coupling \cite{Mehl2016}. For the sake of clarity, we omit both additional variants in this contribution. 
To achieve a stable coupling also for strongly coupled problems and to increase accuracy, we consider an implicit (or strong) coupling, i.e., we repeat $\DSC$ and $\NSC$ multiple times per time window until convergence of $\cDe$. This corresponds to a fixed-point iteration for Eq.~(\ref{eq:FPSC}). If we execute both solvers only once per time step, this corresponds to classical Godunov splitting. It always yields only first-order in time, even if we apply implicit coupling iterating Eq.~(\ref{eq:FPSC}) until convergence. 
We can summarize multi-rate single value coupling in an algorithmic setting as
  \begin{enumerate}
    \item Initialize $\cD^0$ with the value $\cDe$ from the previous time step.
    \item \label{solveDirichlet} Solve the Dirichlet problem (with $\nD$ time steps): $\cN^{k+1} = \DSC( \cD^k )$.
    \item Solve the Neumann problem (with $\nN$ time steps): $\cD^{k+1} = \NSC (\cN^{k+1} )$. 
    \item If $\left\|\cD^{k+1} - \cD^k \right\|$ is small enough: Go to the next time window.
    \item \label{goback:SC} Else: $k=k+1$ and go back to \cref{solveDirichlet}.
 \end{enumerate}
  
Note that, in general, this iteration does not converge and requires an additional acceleration step such as underrelaxation or quasi-Newton in item \ref{goback:SC} in the algortihm above. We describe details for this
acceleration in Sect.~\ref{sec::quasi-newton}.

\subsection{Extension to Waveform Iteration}
\label{sec:waveform}

Waveform iterations (WI) are an elegant way to increase the order with minimal use of the solvers' internal details, i.e., well suited for black-box coupling. In contrast to the method from the previous section, WI exchange more data than just the solutions at the end points. 

The basic idea is to use a time-continuous representation of $\cD$ instead of only the value $\cDe$ at the end of the coupling window to connect the two solvers. We can construct representations of the continuous functions $\cD$ and $\cN$ from all time steps in the window $[t_{\text{ini}}; t_{\text{ini}} + \Delta t]$ by interpolation based on a finite dimensional basis. We, thus, define $\cD$ and $\cN$ as functions in time:
\begin{equation*}
\cD: [t_{\text{ini}}; t_{\text{ini}} + \Delta t] \rightarrow \Gamma \; \mbox{ and } \; \cN: [t_{\text{ini}}; t_{\text{ini}} + \Delta t] \rightarrow \Gamma.
\end{equation*}

If both solvers execute independent time steps within the time window as introduced above, they can now evaluate the continuous representations $c_\mathcal{D}$ and $c_\mathcal{N}$ as boundary conditions wherever necessary. If the Dirichlet solver, for example, executes three time steps using the trapezoidal rule and $\delta t_\mathcal{D} = \frac{1}{3} \Delta t$, it samples $\cD$ at four points in time ($t_{\text{ini}}$, $t_{\text{ini}} + \frac{1}{3} \Delta t$, $t_{\text{ini}} + \frac{2}{3} \, \Delta t$, $t_{\text{ini}} + \Delta t$). If the Neumann solver uses two implicit Euler time steps of size $\delta t_\mathcal{N} = \frac{1}{2} \Delta t$, it evaluates $\cN$ at two points in time ($t_{\text{ini}} + \frac{1}{2} \, \Delta t$ and $t_{\text{ini}} + \Delta t$). If one of the solvers uses higher-order time integration (such as a Runge-Kutta RK4 scheme), necessary substeps can now also be retrieved from $\cD$ or $\cN$. Giving continuous boundary value representations to our solvers yields new operators
\begin{equation*}
\DWI: \cD \mapsto \cN \; \mbox{ and } \; \NWI: \cN \mapsto \cD\;.
\end{equation*}
These new operators are composed of two parts, the actual solver actions
\begin{equation*}
\DS: \cD \mapsto \left( \cN_1, \ldots, \cN_{\nD} \right) \; \mbox{ and } \;
   \NS: \cN \mapsto \left( \cD_1, \ldots, \cD_{\nN} \right)
\end{equation*}
producing time-discrete data at the end of each of their $\nD$ or $\nN$ time steps 
from a time-continuous boundary condition and interpolation operators\footnote{We occasionally drop the subscript for the quantities $\nD, \nN, \cD, \cN$ if we do not refer to a specific side of the coupled problem, but to either.}
\begin{equation*}
\I_{p,c_0}: \left( c_1, \ldots, c_n \right) \mapsto c
\end{equation*}
generating piecewise polynomial interpolation of degree $p$ from $n$ time-discrete points and an 'initial value' $c_0$ at $t_{\text{ini}}$ (i.e., a total of $n' = n + 1$ datapoints). Here, we use an interpolating B-Spline curve with $n'+p+1$ knots \cite{Dierckx1993}\footnote{Note that the first $p$ knots are only artificially added times $< t_{\text{ini}}$ with no associated function values. These additional points in time are used only for the recursive construction of the basis functions at the knots $p+1, \ldots, n'$, our actual interpolation points.}:
\[ \DWI = \I_{p,\cN_0} \circ \DS \; \mbox{ and } \; \NWI = \I_{p, \cD_0} \circ \NS. \]
As we only use $n+1$ data points in time and, in particular, neither use values from previous windows nor derivatives in time for interpolation, we have the restriction $p \le n$ for $\I_{p,c_0}$.
We use piecewise linear ($p=1$), quadratic ($p=2$), or cubic ($p=3$)
interpolation. Piecewise linear interpolation requires
two data points, i.e., $n_\mathcal{D} \ge 1, n_\mathcal{N} \ge 1$, and is, thus, always possible, piecewise quadratic requires at least three samples, i.e., $n_\mathcal{D} \ge 2, n_\mathcal{N} \ge 2$, piecewise cubic at least four samples, i.e., $n_\mathcal{D} \ge 3, n_\mathcal{N} \ge 3$. 

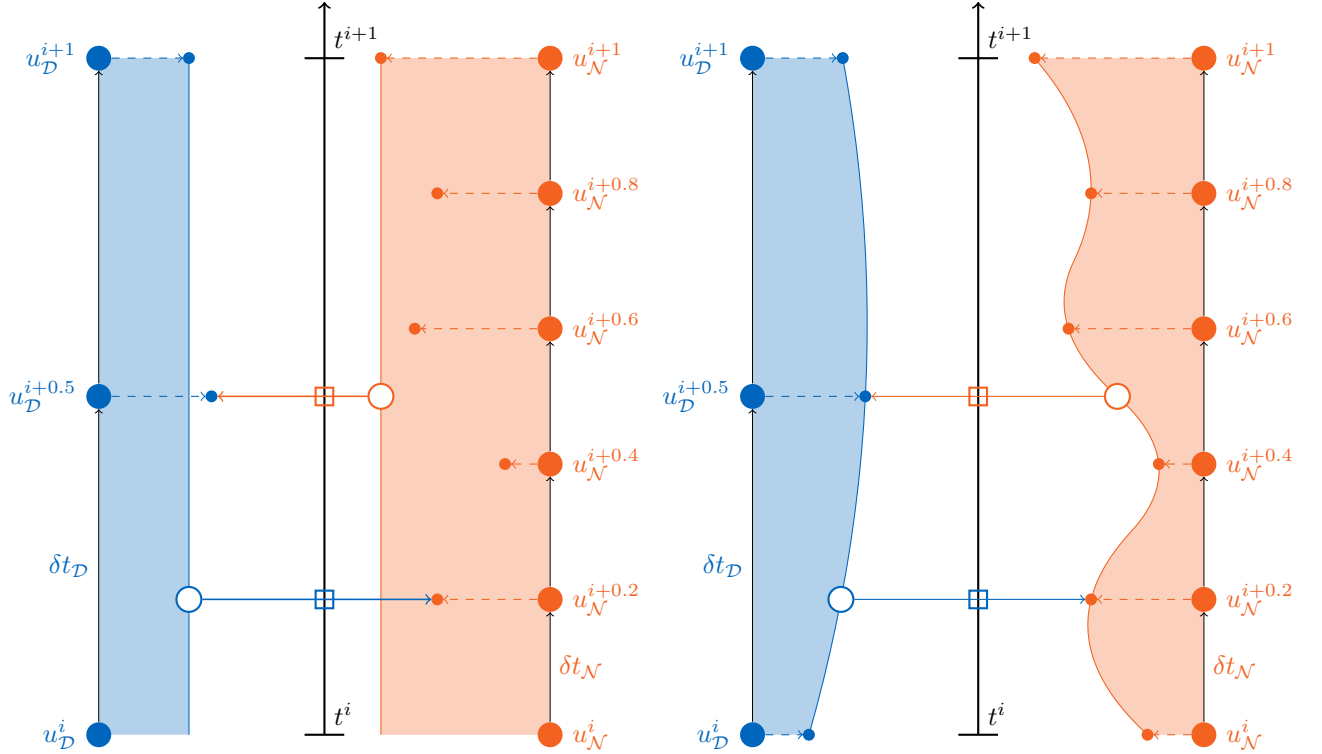
\begin{figure}[h!]
\begin{tikzpicture}[scale=1.5]

\coordinate(t0) at (0,0);
\coordinate(t1) at (0,6);
\coordinate(tsampleA) at (0,1.2);
\coordinate(tsampleB) at (0,3);
\path let \p1 = (t0), \p2 = (t1) in coordinate(t05) at (0,0.5*\y1+0.5*\y2);
\path let \p1 = (t0), \p2 = (t1) in coordinate(t02) at (0,0.8*\y1+0.2*\y2);
\path let \p1 = (t0), \p2 = (t1) in coordinate(t04) at (0,0.6*\y1+0.4*\y2);
\path let \p1 = (t0), \p2 = (t1) in coordinate(t06) at (0,0.4*\y1+0.6*\y2);
\path let \p1 = (t0), \p2 = (t1) in coordinate(t08) at (0,0.2*\y1+0.8*\y2);
\coordinate(participantA) at (-2,0);
\coordinate(participantB) at (2,0);

\draw[->,thick](t0)--($(t1)+(0,.5)$);
\draw[thick] ($(t0)+(5pt,0)$) -- node[above right]{$t^i$} ($(t0)-(5pt,0)$);
\draw[thick] ($(t1)+(5pt,0)$) -- node[above right]{$t^{i+1}$} ($(t1)-(5pt,0)$);

\path let \p1 = (participantA), \p2 = (t0) in node(pA0)[circle,fill=participantBcolor,label={left:\textcolor{participantBcolor}{$u_\mathcal{D}^{i}$}}] at (\x1,\y2){};
\path let \p1 = (participantA), \p2 = (t05) in node(pA05)[circle,fill=participantBcolor,label={left:\textcolor{participantBcolor}{$u_\mathcal{D}^{i+0.5}$}}] at (\x1,\y2){};
\path let \p1 = (participantA), \p2 = (t1) in node(pA1)[circle,fill=participantBcolor,label={left:\textcolor{participantBcolor}{$u_\mathcal{D}^{i+1}$}}] at (\x1,\y2){};

\coordinate(pA0sample) at ($(pA0)+(.5,0)$);
\coordinate(pA05sample) at ($(pA05)+(1,0)$);
\coordinate(pA1sample) at ($(pA1)+(.8,0)$);

\fill[thick, participantBcolor] (pA05sample) circle (1.5pt);
\fill[thick, participantBcolor] (pA1sample) circle (1.5pt);

\draw[participantBcolor, name path = interpA](pA1sample |- pA0sample) -- (pA1sample);

\draw[dashed,participantBcolor,->](pA05) -- ([xshift=-1.5pt]pA05sample);
\draw[dashed,participantBcolor,->](pA1) -- ([xshift=-1.5pt]pA1sample);

\fill[on layer = bg, participantBcolor!30](pA1sample |- pA0sample) -- (pA1sample) -- (pA1.center) -- (pA0.center);

\path let \p1 = (participantB), \p2 = (t0) in node(pB0)[circle,fill=participantAcolor,label={right:\textcolor{participantAcolor}{$u_\mathcal{N}^{i}$}}] at (\x1,\y2){};
\path let \p1 = (participantB), \p2 = (t02) in node(pB02)[circle,fill=participantAcolor,label={right:\textcolor{participantAcolor}{$u_\mathcal{N}^{i+0.2}$}}] at (\x1,\y2){};
\path let \p1 = (participantB), \p2 = (t04) in node(pB04)[circle,fill=participantAcolor,label={right:\textcolor{participantAcolor}{$u_\mathcal{N}^{i+0.4}$}}] at (\x1,\y2){};
\path let \p1 = (participantB), \p2 = (t06) in node(pB06)[circle,fill=participantAcolor,label={right:\textcolor{participantAcolor}{$u_\mathcal{N}^{i+0.6}$}}] at (\x1,\y2){};
\path let \p1 = (participantB), \p2 = (t08) in node(pB08)[circle,fill=participantAcolor,label={right:\textcolor{participantAcolor}{$u_\mathcal{N}^{i+0.8}$}}] at (\x1,\y2){};
\path let \p1 = (participantB), \p2 = (t1) in node(pB1)[circle,fill=participantAcolor,label={right:\textcolor{participantAcolor}{$u_\mathcal{N}^{i+1}$}}] at (\x1,\y2){};

\coordinate(pB0sample) at ($(pB0)-(.5,0)$);
\coordinate(pB02sample) at ($(pB02)-(1,0)$);
\coordinate(pB04sample) at ($(pB04)-(.4,0)$);
\coordinate(pB06sample) at ($(pB06)-(1.2,0)$);
\coordinate(pB08sample) at ($(pB08)-(1,0)$);
\coordinate(pB1sample) at ($(pB1)-(1.5,0)$);

\draw[thin,participantAcolor, name path = interpB](pB1sample |- pB0sample) -- (pB1sample);

\fill[thick, participantAcolor] (pB02sample) circle (1.5pt);
\fill[thick, participantAcolor] (pB04sample) circle (1.5pt);
\fill[thick, participantAcolor] (pB06sample) circle (1.5pt);
\fill[thick, participantAcolor] (pB08sample) circle (1.5pt);
\fill[thick, participantAcolor] (pB1sample) circle (1.5pt);

\draw[thin,dashed,participantAcolor,->](pB02) -- ([xshift=1.5pt]pB02sample);
\draw[thin,dashed,participantAcolor,->](pB04) -- ([xshift=1.5pt]pB04sample);
\draw[thin,dashed,participantAcolor,->](pB06) -- ([xshift=1.5pt]pB06sample);
\draw[thin,dashed,participantAcolor,->](pB08) -- ([xshift=1.5pt]pB08sample);
\draw[thin,dashed,participantAcolor,->](pB1) -- ([xshift=1.5pt]pB1sample);

\fill[participantAcolor!30, on layer = bg](pB1sample |- pB0sample) -- (pB1sample) -- (pB1.center) -- (pB0.center) -- cycle;

\path[name path= connectionAB] let \p1 = (participantA), \p2 = (participantB), \p3 = (tsampleA) in (\x1,\y3) -- (\x2,\y3);

\path[name intersections={of=interpA and connectionAB, by=intersectionA}];

\path[name path= connectionBA] let \p1 = (participantB), \p2 = (participantA), \p3 = (tsampleB) in (\x1,\y3) -- (\x2,\y3);

\path[name intersections={of=interpB and connectionBA, by=intersectionB}];

\node[circle,draw,thick,participantBcolor,fill=white, name intersections={of=interpA and connectionAB,total=\t}](sampleA) at (intersectionA){};

\draw[semithick, ->, participantBcolor, name intersections={of=interpB and connectionAB,total=\t}](sampleA.east) -- ([xshift=-1.5pt]pB02sample);

\node[rectangle,draw=participantBcolor,thick] at (tsampleA){};

\node[circle,draw,thick,participantAcolor,fill=white, name intersections={of=interpB and connectionBA,total=\t}](sampleB) at (intersectionB){};

\draw[semithick, ->, participantAcolor, name intersections={of=interpA and connectionBA,total=\t}](sampleB.west) -- ([xshift=+1.5pt]pA05sample);
	
\node[rectangle,draw=participantAcolor,thick] at (tsampleB){};

\draw[very thin,->](pA0) -- node[left, align=center]{\textcolor{participantBcolor}{$\delta t_\mathcal{D}$}} (pA05);
\draw[very thin,->](pA05) -- node[left, align=center]{} (pA1);

\draw[very thin,->](pB0) -- node[right, align=center]{\textcolor{participantAcolor}{$\delta t_\mathcal{N}$}} (pB02);
\draw[very thin,->](pB02) -- node[right, align=center]{} (pB04);
\draw[very thin,->](pB04) -- node[right, align=center]{} (pB06);
\draw[very thin,->](pB06) -- node[right, align=center]{} (pB08);
\draw[very thin,->](pB08) -- node[right, align=center]{} (pB1);

\end{tikzpicture} %
\begin{tikzpicture}[scale=1.5]

\coordinate(t0) at (0,0);
\coordinate(t1) at (0,6);
\coordinate(tsampleA) at (0, 1.2);
\coordinate(tsampleB) at (0, 3.0);
\path let \p1 = (t0), \p2 = (t1) in coordinate(t05) at (0,0.5*\y1+0.5*\y2);
\path let \p1 = (t0), \p2 = (t1) in coordinate(t02) at (0,0.8*\y1+0.2*\y2);
\path let \p1 = (t0), \p2 = (t1) in coordinate(t04) at (0,0.6*\y1+0.4*\y2);
\path let \p1 = (t0), \p2 = (t1) in coordinate(t06) at (0,0.4*\y1+0.6*\y2);
\path let \p1 = (t0), \p2 = (t1) in coordinate(t08) at (0,0.2*\y1+0.8*\y2);
\coordinate(participantA) at (-2,0);
\coordinate(participantB) at (2,0);

\draw[->,thick](t0)--($(t1)+(0,.5)$);
\draw[thick] ($(t0)+(5pt,0)$) -- node[above right]{$t^i$} ($(t0)-(5pt,0)$);
\draw[thick] ($(t1)+(5pt,0)$) -- node[above right]{$t^{i+1}$} ($(t1)-(5pt,0)$);

\path let \p1 = (participantA), \p2 = (t0) in node(pA0)[circle,fill=participantBcolor,label={left:\textcolor{participantBcolor}{$u_\mathcal{D}^{i}$}}] at (\x1,\y2){};
\path let \p1 = (participantA), \p2 = (t05) in node(pA05)[circle,fill=participantBcolor,label={left:\textcolor{participantBcolor}{$u_\mathcal{D}^{i+0.5}$}}] at (\x1,\y2){};
\path let \p1 = (participantA), \p2 = (t1) in node(pA1)[circle,fill=participantBcolor,label={left:\textcolor{participantBcolor}{$u_\mathcal{D}^{i+1}$}}] at (\x1,\y2){};

\coordinate(pA0sample) at ($(pA0)+(.5,0)$);
\coordinate(pA05sample) at ($(pA05)+(1,0)$);
\coordinate(pA1sample) at ($(pA1)+(.8,0)$);

\fill[thick, participantBcolor] (pA0sample) circle (1.5pt);
\fill[thick, participantBcolor] (pA05sample) circle (1.5pt);
\fill[thick, participantBcolor] (pA1sample) circle (1.5pt);

\draw[participantBcolor, domain = 0:6, variable=\x, name path = interpA] plot ({-7/180 * \x * \x + 51/180 * \x + 90/180 - 2},{\x});

\draw[participantBcolor,dashed,->](pA0) -- ([xshift=-1.5pt]pA0sample);
\draw[participantBcolor,dashed,->](pA05) -- ([xshift=-1.5pt]pA05sample);
\draw[participantBcolor,dashed,->](pA1) -- ([xshift=-1.5pt]pA1sample);

\path let \p1 = (participantB), \p2 = (t0) in node(pB0)[circle,fill=participantAcolor,label={right:\textcolor{participantAcolor}{$u_\mathcal{N}^{i}$}}] at (\x1,\y2){};
\path let \p1 = (participantB), \p2 = (t02) in node(pB02)[circle,fill=participantAcolor,label={right:\textcolor{participantAcolor}{$u_\mathcal{N}^{i+0.2}$}}] at (\x1,\y2){};
\path let \p1 = (participantB), \p2 = (t04) in node(pB04)[circle,fill=participantAcolor,label={right:\textcolor{participantAcolor}{$u_\mathcal{N}^{i+0.4}$}}] at (\x1,\y2){};
\path let \p1 = (participantB), \p2 = (t06) in node(pB06)[circle,fill=participantAcolor,label={right:\textcolor{participantAcolor}{$u_\mathcal{N}^{i+0.6}$}}] at (\x1,\y2){};
\path let \p1 = (participantB), \p2 = (t08) in node(pB08)[circle,fill=participantAcolor,label={right:\textcolor{participantAcolor}{$u_\mathcal{N}^{i+0.8}$}}] at (\x1,\y2){};
\path let \p1 = (participantB), \p2 = (t1) in node(pB1)[circle,fill=participantAcolor,label={right:\textcolor{participantAcolor}{$u_\mathcal{N}^{i+1}$}}] at (\x1,\y2){};

\coordinate(pB0sample) at ($(pB0)-(.5,0)$);
\coordinate(pB02sample) at ($(pB02)-(1,0)$);
\coordinate(pB04sample) at ($(pB04)-(.4,0)$);
\coordinate(pB06sample) at ($(pB06)-(1.2,0)$);
\coordinate(pB08sample) at ($(pB08)-(1,0)$);
\coordinate(pB1sample) at ($(pB1)-(1.5,0)$);

\draw[thin, participantAcolor, name path = interpB, domain = 0:1.8, variable=\x] plot ({0.55827887 * \x * \x - 1.08660131 * \x - 0.5 + 2},{\x}) -- 
plot[domain = 1.8:3, variable=\x] ({-8.52396514e-01 * \x * \x + 3.99183007e+00 * \x + -5.07058824e+00 + 2},{\x}) --
plot[domain = 3:4.2, variable=\x] ({0.6672113289806602 * \x * \x -5.1258169934751665 * \x + 8.605882352951152 + 2},{\x}) --
plot[domain = 4.2:6, variable=\x] ({-0.3730936819222191 * \x * \x + 3.612745098056481 * \x -9.745098039237929 + 2},{\x});

\fill[thick, participantAcolor] (pB0sample) circle (1.5pt);
\fill[thick, participantAcolor] (pB02sample) circle (1.5pt);
\fill[thick, participantAcolor] (pB04sample) circle (1.5pt);
\fill[thick, participantAcolor] (pB06sample) circle (1.5pt);
\fill[thick, participantAcolor] (pB08sample) circle (1.5pt);
\fill[thick, participantAcolor] (pB1sample) circle (1.5pt);

\draw[thin,dashed,participantAcolor,->](pB0) -- ([xshift=1.5pt]pB0sample);
\draw[thin,dashed,participantAcolor,->](pB02) -- ([xshift=1.5pt]pB02sample);
\draw[thin,dashed,participantAcolor,->](pB04) -- ([xshift=1.5pt]pB04sample);
\draw[thin,dashed,participantAcolor,->](pB06) -- ([xshift=1.5pt]pB06sample);
\draw[thin,dashed,participantAcolor,->](pB08) -- ([xshift=1.5pt]pB08sample);
\draw[thin,dashed,participantAcolor,->](pB1) -- ([xshift=1.5pt]pB1sample);

\path[name path= connectionAB] let \p1 = (participantA), \p2 = (participantB), \p3 = (tsampleA) in (\x1,\y3) -- (\x2,\y3);

\path[name intersections={of=interpA and connectionAB, by=intersectionA}];
\path[name intersections={of=interpB and connectionAB, by=intersectionBfromA}];

\path[name path= connectionBA] let \p1 = (participantB), \p2 = (participantA), \p3 = (tsampleB) in (\x1,\y3) -- (\x2,\y3);

\path[name intersections={of=interpB and connectionBA, by=intersectionB}];
\path[name intersections={of=interpA and connectionBA, by=intersectionAfromB}];

\node[circle,draw,thick,participantBcolor,fill=white, name intersections={of=interpA and connectionAB,total=\t}](sampleA) at (intersectionA){};

\draw[->, participantBcolor, name intersections={of=interpB and connectionAB,total=\t}](sampleA.east) -- ([xshift=-1.5pt]intersectionBfromA);
	
\node[thick,draw=participantBcolor] at (tsampleA) {};

\node[circle,draw,thick,participantAcolor,fill=white, name intersections={of=interpB and connectionBA,total=\t}](sampleB) at (intersectionB){};

\draw[->, participantAcolor, name intersections={of=interpA and connectionBA,total=\t}](sampleB.west) -- ([xshift=+1.5pt]intersectionAfromB);
	
\node[rectangle,draw=participantAcolor,thick] at (tsampleB){};

\draw[very thin,->](pA0) -- node[left, align=center]{\textcolor{participantBcolor}{$\delta t_\mathcal{D}$}} (pA05);
\draw[very thin,->](pA05) -- node[left, align=center]{} (pA1);

\draw[very thin,->](pB0) -- node[right, align=center]{\textcolor{participantAcolor}{$\delta t_\mathcal{N}$}} (pB02);
\draw[very thin,->](pB02) -- node[right, align=center]{} (pB04);
\draw[very thin,->](pB04) -- node[right, align=center]{} (pB06);
\draw[very thin,->](pB06) -- node[right, align=center]{} (pB08);
\draw[very thin,->](pB08) -- node[right, align=center]{} (pB1);

\fill[participantBcolor!30,on layer=bg]
plot[domain = 0:6, variable=\x] ({-7/180 * \x * \x + 51/180 * \x + 90/180 - 2},{\x}) --
(pA1.center) -- (pA05.center) -- (pA0.center) -- cycle;
\fill[participantAcolor!30, on layer=bg]
plot[domain = 0:1.8, variable=\x] ({0.55827887 * \x * \x - 1.08660131 * \x - 0.5 + 2},{\x}) -- 
plot[domain = 1.8:3, variable=\x] ({-8.52396514e-01 * \x * \x + 3.99183007e+00 * \x + -5.07058824e+00 + 2},{\x}) --
plot[domain = 3:4.2, variable=\x] ({0.6672113289806602 * \x * \x -5.1258169934751665 * \x + 8.605882352951152 + 2},{\x}) --
plot[domain = 4.2:6, variable=\x] ({-0.3730936819222191 * \x * \x + 3.612745098056481 * \x -9.745098039237929 + 2},{\x}) --
(pB1.center) -- (pB0.center) -- (pB0sample) -- cycle;

\end{tikzpicture} \caption{
  Comparison between single value coupling (\textbf{SC($2,5$)}, left) and waveform iteration (\textbf{WI($2,5;2$)}, right). 
  Exemplary sample values are depicted for $t_{\text{ini}}+\delta t_\mathcal{N}$ and $t_{\text{ini}}+\delta t_\mathcal{D}$. Note that \textbf{SC($2,5$)} is not equivalent to \textbf{WI($2,5;0$)}, i.e., waveform iteration with constant interpolation. This is due to the fact that single value coupling does not use piecewise constant interpolation, but globally constant interpolation over the whole time window.}
\label{fig::coupling}
\end{figure}

The definition of the new operators $\DWI$ and $\NWI$ yields a new fixed-point equation operating on functions in time
\begin{equation}
\cD = \NWI \circ \DWI (\cD) \label{eq:FPWI}
\end{equation}
or, expressed in terms of the time-discrete steps
\begin{equation}
\left( \cD_1, \ldots, \cD_{\nN} \right) = \underbrace{\NS \circ \I_{p,\cN_0} \circ \DS \circ \I_{p,\cD_0}}_{\mbox{$=: H_{\nD, \nN}^{p}$}} \left( \cD_1, \ldots, \cD_{\nN} \right). \label{eq:FPWI_d}
\end{equation}
We denote coupling methods solving Eq.~(\ref{eq:FPWI}) or Eq.~(\ref{eq:FPWI_d}) 
as waveform iteration \textbf{WI($\nD,\nN;p$)} and the corresponding fixed-point operator as $H_{\nD, \nN}^{p}$. Figure \ref{fig::coupling} compares \textbf{SC($2,5$)} and \textbf{WI($2,5;2$)} in a schematic way. 

In an algorithmic setting, we can summarize waveform iteration as
  \begin{enumerate}
    \item From the previous window, get $\cD_0$. Set the initial approximation as constant extrapolation $\cD^0 \equiv \cD_0$\footnote{Other choices for initialization are possible, but are beyond the scope of this paper.}. 
    \item \label{solveDirichlet} Solve the Dirichlet problem (with $\nD$ time steps): 
		      $\left(\cN_1^{k+1}, \ldots, \cN_{\nD}^{k+1} \right) = \DS \left( \cD^k \right)$.
    \item Interpolate to obtain a continuous waveform: $\cN^{k+1} = \I_{p, \cN_0} \left(\cN_1^{k+1}, \ldots, \cN_{\nD}^{k+1} \right)$.
    \item Solve the Neumann problem (with $\nN$ time steps):
		      $\left(\cD_1^{k+1}, \ldots, \cD_{\nN}^{k+1} \right) = \NS \left( \cN^{k+1} \right)$.
    \item Interpolate to obtain a continuous waveform: $\cD^{k+1} = \I_{p, \cD_0} \left(\cD_1^{k+1}, \ldots, \cD_{\nN}^{k+1} \right)$. 
    \item If $\left\|\cD^{k+1}(t) - \cD^k(t) \right\|$ small enough: Go to the next window. 
    \item Else: $k=k+1$ and go back to \cref{solveDirichlet}.
  \end{enumerate}
  
Note that, in the general case, interpolation of different degrees $\pD$ and $\pN$ may be used for the two problems.
For the sake of simplicity, we restrict ourselves to the case where $\pD = \pN = p$.
  
At this point, we would like to add some comments on which combinations of interpolation methods and time stepping schemes have to be used to  (i) reproduce the solution of a problem with known polynomial solution of degree $p$ in time exactly and to (ii) achieve a desired approximation order $p$ in time in general:

For (i), we have to use piecewise polynomial interpolation of at least degree $p$ to represent this solution exactly at the interface. For time integration, a scheme of at least order $p$ has to be used.

For (ii), we have to use a $p$th-order polynomial interpolation, since we can at most afford an error of $\mathcal{O}(\Delta t^{p+1})$ per time window, which accumulates to an error of $\mathcal{O}(\Delta t^p)$ over 
$\mathcal{O}(\Delta t^{-1})$ time steps. For time integration, again a scheme of at least order $p$ has to be used.

Accuracy is in a way limited by the number of time steps $n$ that defines the maximal polynomial degree $p$ of the interpolation scheme. In future work, we intend to remove this restriction by using data from previous windows.

Also this iteration needs acceleration methods that we describe in a general formulation in Sect.~\ref{sec::quasi-newton} and more specifically for waveform iteration in Sect.~\ref{sec::qn-wc}.

\subsection{Interface Quasi-Newton Methods}
\label{sec::quasi-newton}

After having introduced two specific fixed-point operators $H_{\text{SC}}$ and $H_{\nD, \nN}^{p}$ in the previous two sections, we now consider a general fixed-point equation and the associated residual,  
\[H(x)=x, x \in \mathbb{R}^m \;\;\text{and}\;\; R(x) = H(x) - x = 0\;,\] and explain how we can formulate an accelerated iteration 
\[x^{k+1} = \mathcal{A}(H(x^k))\,,\]
with $\mathcal{A}$ denoting the acceleration operator.
We stop the iteration, when the norm of the coupling residual $\|R(x^k)\|$ is small enough.
If $\mathcal{A}$ is the identity, we get a classical fixed-point iteration (Picard iteration), which only converges if $H$ is a contraction. For FSI, this is typically not true \cite{VanBrummelen2009_AddedMassEffect}. In early FSI days, the acceleration operator $\mathcal{A}$ has been realized as underrelaxation,
\[
x^{k+1} = \omega H(x^k) + (1-\omega) x^k\,,\omega \in ]0;1]\,,
\]
e.g., based on a dynamic Aitken scheme \cite{Kuettler2008}. A better approach is to reuse past iterates to build up an approximate Jacobian of $R$ and use this in quasi-Newton iterations \cite{Degroote2009_Performance, Haelterman2009_QNLS}. To avoid linear dependencies between information from previous iterations, we have to perform a modified Newton iteration starting from the result of the pure fixed-point iteration (for details, see 
\cite{Uekermann2016}):
\begin{equation*}
  x^{k+1} = \tilde{x}^k + \Delta \tilde{x}^k \;\;\text{with}\;\; \Delta \tilde{x}^k = J_{\tilde{R}}^{-1}(\tilde{x}^k) \tilde{R}(\tilde{x}^k),
\end{equation*}
where we use $\tilde{x}:=H(x)$, the modified residual $\tilde{R}(\tilde{x}):= \tilde{x} - H^{-1}(\tilde{x})$, and $J_{\tilde{R}}^{-1}(\tilde{x}^k)$ as an approximation of the inverse of the Jacobian of $\tilde{R}$. To compute $J_{\tilde{R}}(\tilde{x}^k)$, we collect input-output information from past iterates\footnote{We use a simple underrelaxation for the first iteration.} of $H$ in tall and skinny matrices $V_k, W_k \in \mathbb{R}^{m \times k}$, $m \gg k$,
\begin{eqnarray*}
   V_k &=& \left[\tilde{R}(\tilde{x}^1) - \tilde{R}(\tilde{x}^0), \tilde{R}(\tilde{x}^2) - \tilde{R}(\tilde{x}^1), \ldots, \tilde{R}(\tilde{x}^k) - \tilde{R}(x^{k-1}) \right]\;, \\
   W_k &=& \left[\tilde{x}^1 - \tilde{x}^0, \tilde{x}^2 - \tilde{x}^1, \ldots, \tilde{x}^k - \tilde{x}^{k-1} \right]\;.
\end{eqnarray*}
For transient coupled problems, we need to solve a fixed-point equation per time window. In this case, it is beneficial for efficiency to also include iterates from past windows in $V_k$ and $W_k$\cite{Uekermann2016}. This variant is however beyond the scope of this paper and subject to future work.
The matrices $V_k$ and $W_k$ define the multi-secant equations for the inverse Jacobian 
\begin{equation}
\label{equ::secant}
J_{\tilde{R}}^{-1} (\tilde{x}^k) \; V_k = W_k\;.
\end{equation}
To get the classical IQN-ILS \cite{Degroote2009_Performance}, we close Eq.~(\ref{equ::secant}) by 
\[\|J^{-1}_{\tilde{R}}(\tilde{x}^k)\|_F \rightarrow \text{min} \,,\]
which gives
\begin{equation*}
J^{-1}_{\tilde{R}} = W_k (V_k^T V_k)^{-1} V_k^T \,.
\end{equation*}
In practice, the Jacobian is not stored, but calculated in matrix-free form from via a QR-decomposition of $V_k$,
\begin{equation}
\label{equ::least-squares}
\alpha = \text{argmin}_{\hat{\alpha} \in \mathbb{R}^k} \|V_k \hat{\alpha} + R^k
\|_2 \;\; \text{followed by} \;\; \Delta \tilde{x}^k = W_k \alpha \;.
\end{equation}

\subsection{Quasi-Newton Waveform Iteration}\label{sec::qn-wc}

Now we combine the two ideas we discussed above: (i) We have used a waveform iteration approach dealing with multi-rate settings, i.e., different numbers $\nD$ and $\nD$ of time steps of size $\delta t_\mathcal{D}$ and $\delta t_\mathcal{N}$ within a window of size $\Delta t$ in the involved solvers, such that $\Delta t = \nN \delta t_\mathcal{N} = \nD \delta t_\mathcal{D}$, and, at the same time, maintaining higher order. (ii) We have presented the concept of quasi-Newton methods as an acceleration approach for the resulting interface fixed-point equations. Now, we show how we can combine the two parts towards a robust, efficient, and accurate coupling for partitioned simulations with black-box solvers.

The missing piece here is the answer to the question which residual components we should consider to build our quasi-Newton scheme. 
For a multi-rate setting with given $\nD$ and $\nN$ based on the fixed-point equation Eq.~(\ref{eq:FPSC}), the choice is canonical, i.e., we use $x^k:= \cDe^k$. %
If we have $m_\mathcal{D}$ degrees of freedom at the interface, $\cDe \in \mathbb{R}^{m_\mathcal{D}}$, we get $V_k, W_k \in \mathbb{R^{m_\mathcal{D}\, \times \, k}}$. We call this variant \textbf{QN-SC}. 

For waveform iteration, we have formulated our Dirichlet-Neumann coupling in terms of the fixed-point equation Eq.~(\ref{eq:FPWI_d}).
We have all $\nD$ time steps of $\cD$ included in the fixed-point equation Eq.~(\ref{eq:FPWI_d}) with the respective fixed point operator $H_{\nD, \nN}^{p}$. The canonical quasi-Newton variant for this system concatenates all time steps' data in a large vector for both input and output of the operator $\tilde{R}$. Thus, we get matrices $V_k$ and $W_k$ in $\mathbb{R}^{(m_\mathcal{D} \cdot \nD) \times k}$. We call this variant \textbf{QN-WI}.
This method is computationally more expensive as we have to deal with a substantially larger interface system than for single value coupling. Typically, the size of the interface system has, however, a negligible influence on the overall performance\cite{Uekermann2016}.

A variant to reduce the size of multi-secant system Eq.~(\ref{equ::secant}) is to measure the residuals not at all samples, but only a subset.
To test such a setting in the following, we only consider one specific case: we only use the last sample $\cD_{\nD}$, i.e.,
\begin{eqnarray*}
  V_k &=& \left[(\tilde{R}(\tilde{x}^1))_{\nD} - (\tilde{R}(\tilde{x}^0))_{\nD}, (\tilde{R}(\tilde{x}^2))_{\nD} - (\tilde{R}(\tilde{x}^1))_{\nD}, \ldots, (\tilde{R}(\tilde{x}^k))_{\nD} - (\tilde{R}(x^{k-1}))_{\nD} \right]\; \in \mathbb{R}^{m_\mathcal{D} \times k}, \\
  W_k &=& \left[\tilde{x}^1 - \tilde{x}^0, \tilde{x}^2 - \tilde{x}^1, \ldots, \tilde{x}^k - \tilde{x}^{k-1} \right]\; \in \mathbb{R}^{(m_\mathcal{D} \cdot \nD) \times k}.
\end{eqnarray*}
In this case, the approximate matrix $J_{\tilde{R}}^{-1}$ is no longer a square matrix but in $\mathbb{R}^{(m_\mathcal{D} \cdot \nD) \times m_\mathcal{D}}$.
However, we do not explicitly calculate $J_{\tilde{R}}^{-1}$, but solve a least squares problem as in Eq.~(\ref{equ::least-squares}) with $\alpha \in \mathbb{R}^k$ as before. We call this variant reduced quasi-Newton \textbf{rQN-WI}.

\textbf{Remark 1:} If we do not include all time steps in the input vector $x$ and, thus, in the matrix $W_k$, but simply restrict the fixed-point iteration formulation entirely to the last step $\nD$, we would also only update the entry $\cD_{\nD}$ and, thus, not generate an appropriate update of the whole waveform polynomial representation of $\cD$.

\textbf{Remark 2:} The method of De Moerlosse et al.\cite{DeMoerloose2019} corresponds to \textbf{WI($1,n;1$)} or \textbf{WI($n,1;1$)}, i.e., one solver executes $n$ time steps $\delta t$ within one window, while the time step of the other solver is equal to the window size $\Delta t$. Linear interpolation is used to transfer information from the larger time step ($\Delta t$) solver to the other solver ($\delta t = \Delta t / n$). In this case, \textbf{QN-WI} and \textbf{rQN-WI} are equivalent.
In our general setting, we impose no restrictions on time steps of the involved solvers except that they 'meet' at a common time at the end of the window, i.e., the window size $\Delta t$ is a multiple of the solver's time steps $\delta t$. Therefore, our generalization also allows to naturally handle more than two coupled solvers.

\textbf{Remark 3:} In literature, waveform iteration, i.e., $H_{\nD, \nN}^{p}$, is usually combined with (constant) underrelaxation. We compare our quasi-Newton variants to this standard technique in Sect.~\ref{sec::results} and call the latter \textbf{rel-WI}. 

Summarizing, we can state that \textbf{QN-WI} and \textbf{rQN-WI} can be considered as quasi-Newton variants in space-time, whereas \textbf{QN-SC}
only acts in space.
\section{Software}\label{sec::software}

We now describe how we realize the aforementioned multi-rate coupling algorithms in software. 
In fact, the bigger picture of our research is to implement an elaborate and sophisticated \textit{treatment of the temporal dimension} into the widespread open-source coupling software preCICE \cite{Bungartz2016_preCICE}. This manuscript is an important milestones for this endeavour.
preCICE is a library that allows to couple existing (legacy) simulation codes to complete multiphysics simulations in a minimally invasive way. preCICE is written in \texttt{C++}, but also offers bindings for \texttt{C}, \texttt{Fortran} and \texttt{Python}. 
For many community codes (such as OpenFOAM, CalculiX or SU2), ready-to-use adapters are provided \cite{Uekermann2017_Adapters}. 
The setup of coupled codes and of the data to be exchanged is configured at run-time through a global XML file.
Various quasi-Newton algorithms are already implemented in preCICE including IQN-ILS \cite{Degroote2009_Performance} as introduced in Sect.~\ref{sec::quasi-newton}. 

Whereas single value coupling, as introduced in Sect.~\ref{sec:basic}, is based on exchanging the values of the coupling variables at window end points, waveform iteration, as introduced in Sect.~\ref{sec:waveform} needs time-continuous representation of the coupling variables over complete windows. This important mathematical difference is also reflected in the software requirements. 
preCICE currently only offers interface methods such as \texttt{data = read\_data(...)} and \texttt{write\_data(data,...)} to read and write values at a single point in time, which can be interpreted in most cases as the window end $t_{\text{ini}}+\Delta t$. 
Afterwards, the interface method \texttt{advance()} triggers the exchange of these values between both solvers and the computation of the quasi-Newton acceleration. This set of methods provides everything needed to realize single value coupling, which, in fact, is already fully supported\cite{Gatzhammer2015} by the preCICE version we use in this contribution, v1.6.1 . Waveform iteration, however, requires some notion of time in the application programming interface (API) of preCICE, which would make a fundamental redesign of the coupling logic necessary.

\paragraph{Waveform Bindings}

\begin{figure}
\centering
\begin{tikzpicture}[scale=2]

\draw[draw=none] (2.5,-.7) rectangle (10.5,1.7);

\coordinate(fenicsOrigin) at (6,0);
\draw[fill=participantAcolor](fenicsOrigin) rectangle node{Solver} ++(2.3, 1);

\coordinate(adapterOrigin) at (6,0);
\draw[fill=participantAcolor](adapterOrigin) -- node(adapterA)[below]{} ++ (-.75,0) -- ++(0,.25) -- ++ (.375,.25) -- ++(-.375,.25) -- ++(0,.25) -- ++(.75,0) -- cycle;

\coordinate(preciceOrigin) at (3.85,0);
\draw[fill=preciceColor](preciceOrigin) -- ++(.34,0) node[below,pos=.5](preciceaim){}  -- ++(0,.27) -- ++(0.375,0) -- ++(0,.46) -- ++(-.375,0) -- ++(0,.27) -- ++(-.34,0) -- cycle;

\coordinate(WRlayerOrigin) at (4.55,0);
\draw[fill=red!50](WRlayerOrigin) -- ++(.5,0) -- ++(0,.25) -- ++ (.375,.25) -- ++(-.375,.25) -- ++(0,.25) -- node[above,pos=.5](wraim){} ++ (-.75,0) -- ++(0,-.27) -- ++(0.375,0) -- ++(0,-.46) -- ++(-.375,0) -- ++(0,-.27) -- ++(.75,0) -- cycle;

\node(precice) at (4.25,-.7){\texttt{libprecice}: \texttt{u = read\_vector\_data("Forces1")}};
\draw[](precice) -- (preciceaim);

\node(wrprecice) at (6,1.5){\texttt{waveform-bindings}: \texttt{u = read\_vector\_data("Forces", time)}};
\draw[](wrprecice) -- (wraim);

\node[align=left](fenicsadapter) at (6,-.5){\texttt{solver adapter}};
\draw (fenicsadapter) -- ($(adapterOrigin)-(0.25,0.1)$);

\end{tikzpicture} \caption{Sketch of layered software architecture for waveform iteration. Instead of directly calling the preCICE API from the solver adapter, the waveform bindings are called which in turn call preCICE in addition to executing new components such as interpolation in time. Data access methods need to be augmented by time as input argument.}
\label{fig:waveformbindings}
\end{figure}
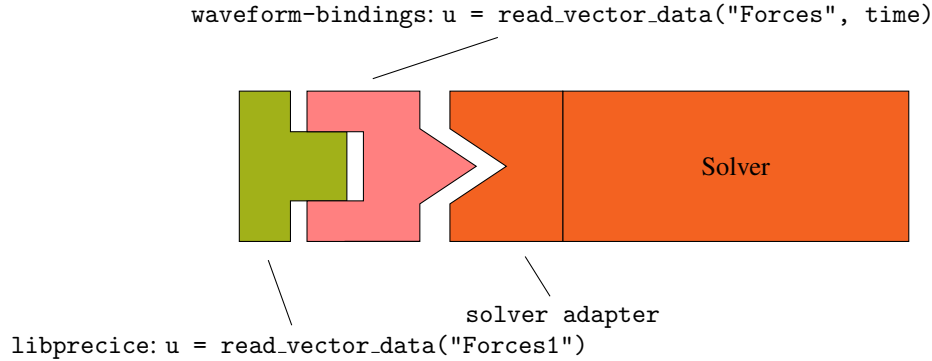

As the choice of the most suitable waveform iteration algorithm for preCICE is still an open research question, we instead aim for a prototype software design, which allows to easily test various multi-rate approaches without requiring immature refactoring of preCICE. A mature implementation of waveform iteration in preCICE is subject to future work. We, therefore, decided that the prototype should implement waveform iteration in a new middle layer called \textit{waveform bindings}, sitting between the adapter code and the preCICE library. Whereas preCICE does the communication and the quasi-Newton acceleration of the coupling data, the waveform bindings handle the time-continuous representation. 
With this solution, we do not require any alterations to preCICE itself. Figure \ref{fig:waveformbindings} sketches the layered architecture.

The essential trick to realize this separation of concerns is to treat a complete time window like a single time step in preCICE. 
To this end, the data fields of preCICE need to be altered: instead of one data field, such as \texttt{Displacements}, we now have one data field per time step of a window, \texttt{Displacements1}, \texttt{Displacements2}, and so forth. preCICE can be configured to handle and communicate multiple fields, but has no notion of time for these fields, they are just several data fields for a single preCICE time step. 
The waveform bindings, on the other hand, have this notion of time: \texttt{Displacements1} are the values after the first time step within a window, \texttt{Displacements2} after the second time step, and so forth.
The API of the waveform bindings now offers methods such as \texttt{write\_data(data,..., time)} to write timestamped data at a specific time step $t \in \{t_{\text{ini}}+ \delta t, t_{\text{ini}}+ 2 \delta t, \ldots, t_{\text{ini}} + \Delta t\}$ and   
\texttt{data = read\_data(..., time)} to read data at an arbitrary $t \in [t_{\text{ini}}, t_{\text{ini}}+\Delta t]$. If the latter does not coincide with a given time step, an interpolated value is returned. We use the SciPy package\footnote{SciPy version 1.3.1} to compute the interpolant, a B-Spline of degree one to four using the SciPy function \texttt{interpolate.splrep}. 
The waveform bindings also offer a function \texttt{advance()} which in turn calls the preCICE \texttt{advance()} upon completion of a window to communicate the discrete values (such as \texttt{Displacements1}, \texttt{Displacements2}, \ldots) between both solvers. In the preCICE configuration, the quasi-Newton acceleration can now be configured to be computed on the full set of data fields or only a reduced subset as explained in Sect.~\ref{sec::qn-wc}.

\paragraph{Application Programming Interface}

To only require minimal modification to existing adapters and testcases, the main API design goal of the waveform bindings\footnote{Code version used: \url{https://github.com/BenjaminRueth/waveform-bindings/commit/3dd475176b5b3540b3b91766957ef0394e816ec1}} is to mimic the preCICE API as closely as possible. In fact, solely the data access functions are altered as explained above to include time as an additional input argument. Therefore, rewriting an existing coupled code to use the waveform bindings instead of preCICE directly requires modifying less than ten lines of code, which makes alongside testing of single value coupling and waveform iteration easy. 

Listing \ref{code:example} exemplarily shows how the waveform bindings are used to couple a structure solver. The code is nearly identical to a coupled code directly using the preCICE API. In the following, we explain the slight differences. For a detailed explanation of the preCICE API itself, we refer to the respective documentation\footnote{\url{https://github.com/precice/precice/wiki/Adapter-Example}}. 
The obvious change in line 1 is to import \texttt{waveform\_bindings} instead of \texttt{precice}. Afterwards, the waveform bindings need additional configuration in lines 11 to 13. The number of time steps per window on both sides as well as the interpolation strategy need to be defined. 
For the data access in lines 23 and 31, the preCICE API is extended as mentioned previously. 
In the example, the structure solver only samples force values at the end of each time step $t + \delta t$ in line 23, as we use an implicit Newmark-beta method and, thus, only require these values. For other time stepping methods, forces could be read from the waveform bindings wherever necessary. 

\begin{listing}
\caption{A simple structure solver in Python using the waveform bindings for \textbf{WI($2,5;1$)}}
\label{code:example}
\begin{minted}[mathescape,linenos,numbersep=5pt,gobble=0,frame=none,framesep=20mm,escapeinside=||,breaklines]{python}
import waveform_bindings as precice

solver_name = "StructureSolver"  # Structure solver is the Neumann solver $\mathcal{N}$
|{\color{orange}interface}| = precice.Interface(solver_name)  # create handle to waveform bindings API
|{\color{orange}interface}|.configure("precice-config.xml")  # usual preCICE configuration

mesh_ID, vertex_IDs = ...  # define coupling mesh on $\Gamma$
u = setup_solver()  # initialize displacement field $u_\mathcal{N}$

# configuration for computing WI(2,5;1)
n_this, n_other = 5, 2  # define $\nN$ and $\nD$
interpolation_strategy = "linear"  # use piecewise linear interpolation with $p=1$
|{\color{orange}interface}|.configure_waveform_iteration(n_this, n_other, interpolation_strategy)

window_size = |{\color{orange}interface}|.initialize()  # preCICE couples at window size $\Delta t$
dt = window_size / n_this  # Neumann solver $\mathcal{N}$ performs time step of size $\delta t_\mathcal{N} = \Delta t / \nN$

while |{\color{orange}interface}|.is_coupling_ongoing():
    if |{\color{orange}interface}|.writing_checkpoint_is_required():  
        write_checkpoint(t, u)  # at window start $t_\text{ini}$

    # read $c_\mathcal{N}$ at $t + \delta t_\mathcal{N}$ needed for this time step from coupling interface
    forces = |{\color{orange}interface}|.read_block_vector_data("Forces", mesh_ID, vertex_IDs, t + dt)
    coupling_boundary_condition = set_forces(forces)  # apply $c_\mathcal{N}$

    # perform time step of size $\delta t_\mathcal{N}$
    u_new = do_time step(u, coupling_boundary_condition, t, dt)
    
    # write values of $c_\mathcal{D}$ at $t+ \delta t_\mathcal{N}$ to coupling interface
    displacements = compute_displacement(new_u)  # extract $c_\mathcal{D}$
    |{\color{orange}interface}|.write_block_vector_data("Displacements", mesh_ID, vertex_IDs, displacements, t + dt)
    
    # advance coupling interface by $\delta t_\mathcal{N}$
    |{\color{orange}interface}|.advance(dt)

    if |{\color{orange}interface}|.reading_checkpoint_is_required():  # rollback to $t_\text{ini}$
        t, u = read_checkpoint()
    else:  # continue to next time step $t + \delta t_\mathcal{N}$
        t, u = t + dt, u_new

|{\color{orange}interface}|.finalize()

\end{minted}
\end{listing}

\paragraph{Limitations}

The layered architecture, however, also comes with technical shortcomings. Firstly, as the waveform bindings build upon the \texttt{Python} bindings of preCICE, currently only codes in \texttt{Python} can be coupled.  
Besides this language restriction, the software concept comes with a few additional limitations as expected for a prototype. For quasi-Newton acceleration, preCICE allows to reuse iterations from previous time steps to achieve faster convergence. As preCICE has no notion of time for the samples \texttt{Displacements1}, \texttt{Displacements2}, \ldots, the reuse strategy would give wrong values. We, therefore, simply switch off this functionality for all experiments in the next section. For the same reason, we cannot use underrelaxation in preCICE, neither directly nor as starting value for the quasi-Newton methods. Otherwise the first relaxed iteration would be a weighted combination between the current solution and the waveform of the last window. Waveform relaxation \textbf{rel-WI} as introduced in Sect.~\ref{sec::qn-wc}, however, needs a weighted combination between the current solution and a constant waveform extrapolated from the end value of the previous window. As a remedy, we realized underrelaxation directly in the \texttt{advance()} function of waveform bindings and switched it off in preCICE. Last, each coupled solver needs to know the number of time steps per window for \textit{both} sides at configuration, see Listing \ref{code:example}, lines 11-13, contradicting the flexibility argument of black-box coupling. We aim to eliminate all these limitation with a mature implementation of waveform iteration within preCICE in future work. Within preCICE in particular means to maintaining backwards compatibility. To which extent this is possible is an open research question.

\section{Numerical Results}
\label{sec::results}

We consider two test cases to evaluate the performance of the proposed algorithms both in terms of convergence properties of the quasi-Newton iterations and in terms of the
achieved approximation order in time: (i) a simple 2D heat equation (diffusion equation) with
known analytic solution, where the computational domain has been divided
artificially into two partitions, (ii) a fluid-structure interaction scenario, where a solid truss is deformed by the surrounding fluid flow.
Both examples are discretized in space using finite elements. To comply with the restriction of the waveform bindings to Python, we use the finite element frameworks FEniCS\footnote{FEniCS version 2019.1.0} \cite{Logg2012, alnaes2015fenics} and Nutils\footnote{Nutils version \url{https://github.com/evalf/nutils/commit/61fd19ce275dff154b9d3481ccc0d849244544f4}} \cite{Nutils}. We choose two different codes to show-case the black-box methodology, for which testcases and adapters were already available. 
For interfacing between FEniCS and preCICE, we make use of the FEniCS adapter\footnote{Source code at \url{https://github.com/precice/fenics-adapter}. For directly accessing preCICE, we use version \url{https://github.com/precice/fenics-adapter/commit/1ccab476e3923b34daf274e365ad777f9411569f}; for using the waveform bindings, we use version \url{https://github.com/precice/fenics-adapter/commit/61c8e2a613f47d79d03eb786f7446f691d9da2bb}}, which converts between FEniCS and preCICE data structures and realizes the rewinding of time steps for implicit (or strong) coupling. Nutils, on the other hand, does not use an adapter, but directly accesses preCICE or the waveform bindings.
Code and instructions for running the examples can be found on \url{https://gitlab.lrz.de/precice/ijnme2019-experiments}.

\subsection{Conjugate Heat Transfer}
\label{ssec:heat}

With this academic test case, we want to test two things: (i) How well does QN converge for WI, which combination works best? (ii) Can we recover higher order in time? We look at both in this
order.
As most of our conclusions are based on a known exact solution, the numerical methods
in the solvers have been chosen accordingly. In particular, the solvers use
finite elements of sufficiently high order to be able to solve the problem without spatial discretization error. 

\paragraph{Scenario Setup}

We solve the dimensionless 2D time-dependent heat equation
\begin{equation}
  \pder{u}{t} = \Delta u + f\;  \mbox{ on }  \; \Omega = \left[0,2\right] \times \left[0,1\right]\; \mbox{ for time } \; t \in \left[0; T\right] \label{eq:heat}\\
  u=g\; \mbox{ on } \; \partial \Omega \; \mbox{ for time } \; t \in \left[0; T\right] \nonumber
\end{equation}
as described in the FEniCS tutorials book\cite{Langtangen2016}. 
In the following, we extend this example, such that it can be used to assess the quality of the previously discussed algorithms. For details on finite element discretization and solution strategies in FEniCS, we refer the interested reader to the FEniCS tutorials book\cite{Langtangen2016}.

The right-hand side $f$ as well as the boundary conditions in Eq.~(\ref{eq:heat}) are manufactured such that we get the analytical solution
\[ u_\text{exact}(x,y,t) = 1+g(t) x^2+ 3 y^2+1.2 t. \]
with an arbitrary continuous function $g$. We consider the following choices for $g(t)$:
\begin{itemize}
  \item polynomial: $g_{\text{pol}}(t)=(1+t)^\gpolorder$,
  \item trigonometric: $g_{\text{tri}}(t)=\sin(t)$.
\end{itemize}

Using $P^2$ finite elements and a sufficient approximation order in time
allows us to recover $u_\text{exact}$ and to compute the exact heat flux
$\vec{q}_\text{exact} = \vec{\nabla} u_\text{exact}$ via finite element projection.

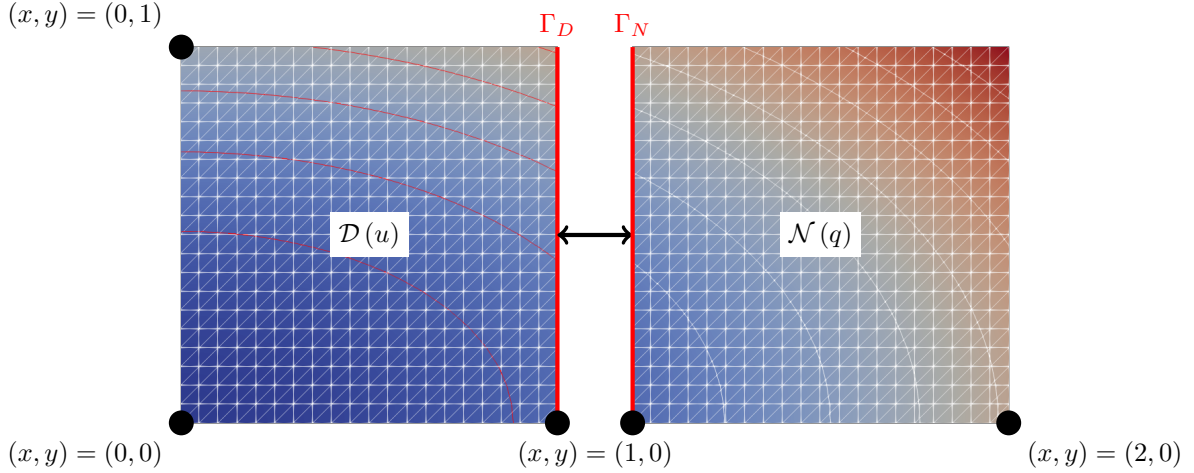
\begin{figure}[h!]
\centering
\begin{tikzpicture}[scale=5]

\draw[fill=participantAcolor!50](0,0) rectangle node{\includegraphics[width=5cm]{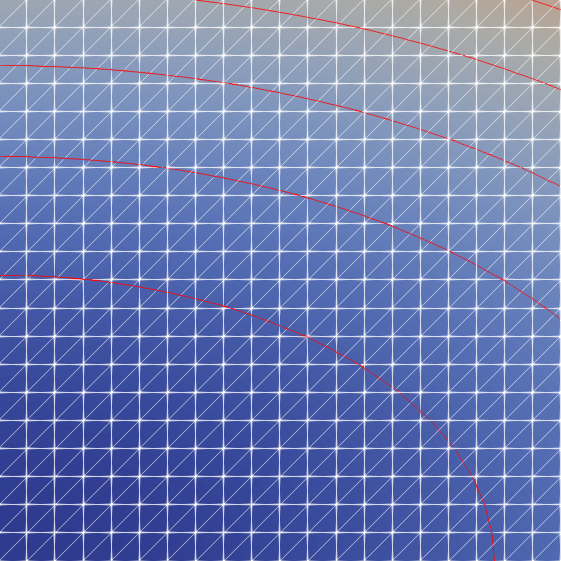}} node[fill=white]{$\mathcal{D}\left(u\right)$}++(1,1);
\draw[ultra thick, red](1,0) -- ++(0,1) node[above]{$\Gamma_D$};

\draw[<->,ultra thick](1,.5) -- (1.2,.5);

\draw[fill=participantBcolor!50](1.2,0) rectangle node{\includegraphics[width=5cm]{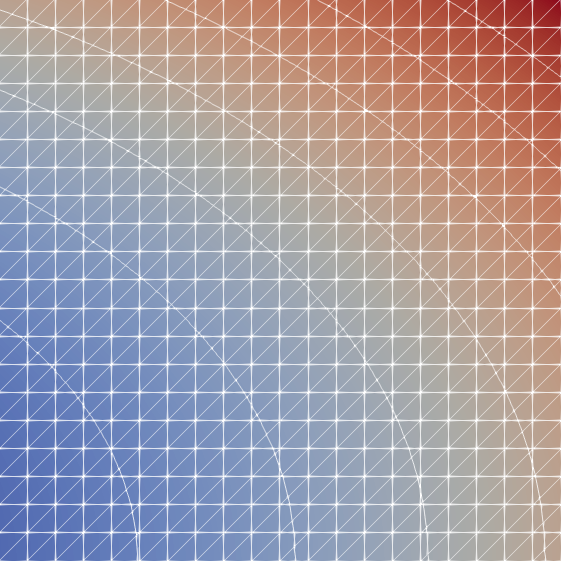}} node[fill=white]{$\mathcal{N}\left(q\right)$} ++(1,1);
\draw[ultra thick, red](1.2,0) -- ++(0,1) node[above]{$\Gamma_N$};
\node[fill, circle, label=below right: {$(x,y)=(2,0)$}] at (2.2,0){};
\node[fill, circle, label=below left: {$(x,y)=(0,0)$}] at (0,0){};
\node[fill, circle, label=above left: {$(x,y)=(0,1)$}] at (0,1){};
\node[label=below: {$(x,y)=(1,0)$}] at (1.1,0){};
\node[fill, circle] at (1,0){};
\node[fill, circle] at (1.2,0){};

\end{tikzpicture} \caption{Heat-heat coupling: partitioning and mesh. Contour lines visualizing temperature are continuous across the coupling interface. This indicates that temperature and heat flux are continuous at the coupling interface.}
\label{fig:heatsetup}
\end{figure}

To use the heat equation as a benchmark for partitioned solution strategies, we divide the domain $\Omega$ into two independent subdomains $\Omega_\mathcal{D} = \left[0,1\right] \times \left[0,1\right]$ (handled by the Dirichlet solver) and $\Omega_\mathcal{N} =  \left[1,2\right] \times \left[0,1\right]$ (handled by the Neumann solver) such that $\Omega = \Omega_\mathcal{D}\cup \Omega_\mathcal{N}$ and the coupling boundary $\Gamma = \Omega_\mathcal{D}\cap \Omega_\mathcal{N}$, compare Fig.~\ref{fig:heatsetup}. This setup will be referred to as the \emph{partitioned heat equation}. We use a triangular equidistant mesh with $20 \times 20 \times2$ triangles in $\Omega_{\mathcal{D}}$ and $\Omega_{\mathcal{N}}$.
We extract the temperature $u$ or the normal component of the heat flux $q_\bot$ on $\Gamma$ from one participant to be provided for the other as Dirichlet boundary condition $u = \cD$ or Neumann boundary condition $\vec{n} \cdot \vec{\nabla} u = q_\bot = \cN$, respectively.

For time stepping, we use implicit Euler (IE), the trapezoidal rule (TR) or a fourth order spectral deferred correction (SDC) scheme \cite{dutt2000spectral}, depending on the choice of $\gpolorder$ and the goal in terms of approximation order in time. We construct the SDC scheme from an implicit Euler time integration scheme with three Gauss-Lobatto collocation nodes at $t_{\text{ini}} + i\delta t, t_{\text{ini}} + (i + 0.5) \delta t$ and $t_{\text{ini}} + (i+1) \delta t$, $i=1\ldots n$ and $K$ correction sweeps per time step of size $\delta t$.

\paragraph{Quasi-Newton Convergence}

In Section \ref{sec::qn-wc}, we introduced two variants combining waveform iteration coupling with quasi-Newton acceleration, %
\textbf{QN-WI} and \textbf{rQN-WI}.
In addition, we consider classical underrelaxation for waveform iteration \textbf{rel-WI} and
the plain fixed-point iteration \textbf{full-WI}.
For \textbf{rel-WI}, we apply an underrelaxation of $0.5$, which is optimal for the
semidiscrete case \cite{gander2016dirichlet}. We compare these methods concerning the speed of
convergence. As a reference, we compare the convergence of all WI coupling acceleration methods to the convergence of \textbf{QN-SC} for \textbf{SC($\nD$,$\nN$)}. The latter quantifies the price we pay for higher order in time
in terms of iteration counts per time step.

As we are not considering the approximation order of our time stepping scheme in this paragraph,
we use implicit Euler time integration combined with piecewise linear interpolation for WI in both
domains (\textbf{WI($\nD$,$\nN;1$)}). Note that for the setups in which waveform iteration is used
(\textbf{QN-WI}, \textbf{rQN-WI}, \textbf{rel-WI}, \textbf{full-WI}), the termination criterion considers interface values at
every time step. For the basic variant with plain single value coupling (\textbf{QN-SC}),
only the values at the end of the time window are considered. We use a moderate relative coupling tolerance of $10^{-5}$ to reach reasonable runtimes.
For all QN acceleration schemes, we use a \textit{QR2 filter}\cite{Haelterman2015_Filtering} for the QR decomposition with limit $\epsilon=10^{-3}$ and a \textit{residual-sum weighting} of the individual sub vectors\cite{Uekermann2016}.

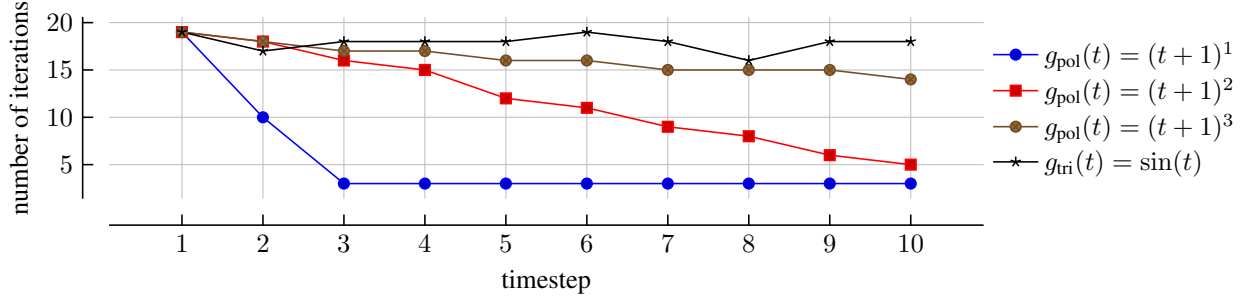
\begin{figure}[h!]
\centering
\begin{tikzpicture}
\begin{axis}[
width=.8\textwidth,
height=4cm,
tuftelike,
legend cell align={left},
xlabel={timestep}, 
ylabel={number of iterations}, 
align=center,
legend style={at={(1,0.5)},anchor=west, draw=none, align=left},
grid
]
\addplot table[x=Timesteps,y=Iterations]{data/residualOverTime/plot_ie/precice-HeatNeumann-iterations_l.log};
\addlegendentry{$g_\text{pol}(t) = (t + 1)^1$};
\addplot table[x=Timesteps,y=Iterations]{data/residualOverTime/plot_ie/precice-HeatNeumann-iterations_q.log};
\addlegendentry{$g_\text{pol}(t) = (t + 1)^2$};
\addplot table[x=Timesteps,y=Iterations]{data/residualOverTime/plot_ie/precice-HeatNeumann-iterations_c.log};
\addlegendentry{$g_\text{pol}(t) = (t + 1)^3$};
\addplot table[x=Timesteps,y=Iterations]{data/residualOverTime/plot_ie/precice-HeatNeumann-iterations_s.log};
\addlegendentry{$g_\text{tri}(t) = \sin(t)$};
\end{axis}
\end{tikzpicture} \caption{Heat-heat coupling: Number of coupling iterations per time step for \textbf{QN-WI}, \textbf{WI($1,1;1$)}, implicit Euler for time stepping with time step/window sizes $\delta t_1 = \delta t_2 = \Delta t = 1$ for $T=10$, different choices for $g(t)$, and a coupling tolerance of $10^{-12}$.}
\label{fig:resovertime}
\end{figure}

As test case, we use $g_{\text{tri}}(t)$, since for $g_{\text{pol}}(t)$, the $g(t)x^2$ term quickly dominates the equation leading to a much simpler coupled problem. See Figure \ref{fig:resovertime} for the iteration per time step for \textbf{QN-WI} for different choices of $g(t)$. Only the periodic $g_{\text{tri}}(t)$ shows a rather constant behavior of the iterations per time step.

We measure the average number of iterations required to terminate in every window until $T=10$ and present the
respective results in Tab.~\ref{tab::qn-iterations}. We evaluate the performance for different window sizes $\Delta t= 0.1, 0.2, 0.5, 1.0, 2.0, 5.0$ and different multi-rate settings with $\nD, \nN \in \{1,3,5\}$. Our study showed that \textbf{full-WI}, i.e., waveform iteration without acceleration, does not converge within 100 iterations per time step. Thus, we present only iteration counts for the other four variants. The most important findings are summarized in the following:

\begin{enumerate}
\item For \textbf{rel-WI}, the number of coupling iterations is almost independent of the multi-rate setup.
      All quasi-Newton variants need two to three times fewer iterations than \textbf{rel-WI}.
\item Overall, \textbf{QN-WI} is more robust, while \textbf{rQN-WI} is more efficient for special cases. In general, \textbf{rQN-WI} seems to have difficulties for setups with many substeps on both sides, while \textbf{QN-WI} can better cope with these setups. \textbf{rQN-WI} and \textbf{QN-WI} perform equally well for \textbf{WI($\nD,1;1$)} as both methods are equivalent for these cases. 
\item For all variants, the number of iterations slowly decreases with decreasing window size $\Delta t$. This effect is observable, but not very significant. Thus, there is potential to make coupling more asynchronous by applying waveform iteration on larger windows.
\item In comparison to existing literature \cite{Uekermann2016}, the number of iterations needed in our experiments is relatively high. We explain the high iteration count with the fact that our algorithms currently do not reuse data of previous time windows. Reusing iterations is a state-of-the-art approach for reducing the number of iterations needed, which remains as future work.
\item \textbf{QN-SC} performs well and does not suffer from multi-rate, even though the time steps ($\delta t$) are ignored by the acceleration scheme. Note, that \textbf{QN-SC} uses an easier convergence criterion than the other two quasi-Newton variants, since only the last sample is considered at the end of the window and recall that this method does not allow to achieve more than first
  order accuracy in time. 
\end{enumerate}

From these findings, we conclude that the convergence speed of the quasi-Newton based variants \textbf{QN-SC}, \textbf{QN-WI} and \textbf{rQN-SC} is relatively similar. The underrelaxation based variant \textbf{rel-WI} performs significantly worse. In the remainder of this paper, we will therefore focus on the quasi-Newton based variants.

\begin{table}[h!]
\centering
\begin{minipage}[t]{0.47\textwidth}
\scalebox{0.75}{
\begin{tabular}{p{2.1cm}|rrrrrr}
\toprule
 \textbf{rel-WI} / $\Delta t$  & \hspace{0.27cm} 5.0 & \hspace{0.27cm} 2.0 & \hspace{0.27cm} 1.0 & \hspace{0.27cm} 0.5 & \hspace{0.27cm} 0.2 & \hspace{0.27cm} 0.1 \\
\midrule
\textbf{WI($1,1;1$)} & \textbf{24.50} &         22.80  & \textbf{21.90} & \textbf{20.85} & \textbf{19.54} & \textbf{18.47} \\
\textbf{WI($1,3;1$)} & \textbf{24.50} & \textbf{22.60} & \textbf{21.90} & \textbf{20.85} & \textbf{19.54} &         18.50 \\
\textbf{WI($1,5;1$)} & \textbf{24.50} &         22.60  & \textbf{21.90} & \textbf{20.85} & \textbf{19.54} &         18.49 \\
\textbf{WI($3,1;1$)} & \textbf{24.50} &         23.00  &         22.10  &         21.05  &         19.64  &         18.51 \\
\textbf{WI($3,3;1$)} &         25.00  &         24.00  &         22.50  &         21.10  &         19.68  &         18.52 \\
\textbf{WI($3,5;1$)} &         25.00  &         24.00  &         22.50  &         21.10  &         19.66  &         18.54 \\
\textbf{WI($5,1;1$)} & \textbf{24.50} &         23.00  &         22.30  &         21.10  &         19.68  & \textbf{18.57} \\
\textbf{WI($5,3;1$)} & \textbf{24.50} & \textbf{24.40} &         22.70  &         21.20  &         19.72  &         18.56 \\
\textbf{WI($5,5;1$)} & \textbf{25.50} & \textbf{24.40} & \textbf{22.80} & \textbf{21.25} & \textbf{19.74} &         18.56 \\
\bottomrule
\end{tabular}
}
\end{minipage}
\hspace{0.02\textwidth}
\begin{minipage}[t]{0.47\textwidth}
\scalebox{0.75}{
\begin{tabular}{p{2.1cm}|rrrrrr}
\toprule
\textbf{QN-SC} / $\Delta t$   & \hspace{0.27cm} 5.0 & \hspace{0.27cm} 2.0 & \hspace{0.27cm} 1.0 & \hspace{0.27cm} 0.5 & \hspace{0.27cm} 0.2 & \hspace{0.27cm} 0.1 \\
\midrule
\textbf{SC($1,1$)} & \textbf{10.50} &         10.00  &         9.00  &         7.85  &         6.54  & \textbf{5.45} \\
\textbf{SC($1,3$)} & \textbf{10.50} & \textbf{10.20} & \textbf{8.60} &         7.75  & \textbf{6.50} &         5.72 \\
\textbf{SC($1,5$)} & \textbf{10.50} &          9.80  & \textbf{8.60} &         7.85  &         6.64  &         5.70 \\
\textbf{SC($3,1$)} & \textbf{10.50} &          9.60  & \textbf{9.30} &         8.10  & \textbf{6.94} &         6.14 \\
\textbf{SC($3,3$)} & \textbf{10.50} &  \textbf{9.20} &         8.80  &         7.65  &         6.58  &         5.88 \\
\textbf{SC($3,5$)} & \textbf{10.50} &  \textbf{9.20} &         8.70  &         7.55  & \textbf{6.50} &         5.87 \\
\textbf{SC($5,1$)} & \textbf{10.50} &          9.60  & \textbf{9.30} & \textbf{8.15} & \textbf{6.94} & \textbf{6.15} \\
\textbf{SC($5,3$)} & \textbf{10.50} &  \textbf{9.20} &         8.80  &         7.60  &         6.62  &         5.94 \\
\textbf{SC($5,5$)} & \textbf{10.50} &  \textbf{9.20} &         8.70  & \textbf{7.45} &         6.52  &         5.92 \\
\bottomrule
\end{tabular}
}
\vspace{0.5cm}
\end{minipage}

\begin{minipage}[t]{0.47\textwidth}
\scalebox{0.75}{
\begin{tabular}{p{2.1cm}|rrrrrr}
\toprule
 \textbf{rQN-WI} / $\Delta t$   & \hspace{0.27cm} 5.0 & \hspace{0.27cm} 2.0 & \hspace{0.27cm} 1.0 & \hspace{0.27cm} 0.5 & \hspace{0.27cm} 0.2 & \hspace{0.27cm} 0.1 \\
\midrule
\textbf{WI($1,1;1$)} & \textbf{10.50} &         10.00  &          9.00  &          7.85  &         6.54  &         5.45 \\
\textbf{WI($1,3;1$)} & \textbf{10.50} &          9.60  &          8.70  &          7.45  &         6.04  &         5.12 \\
\textbf{WI($1,5;1$)} & \textbf{10.50} &  \textbf{9.20} &  \textbf{8.50} &  \textbf{7.30} & \textbf{5.96} & \textbf{4.99} \\
\textbf{WI($3,1;1$)} & \textbf{10.50} &         10.00  &          9.20  &          8.10  &         6.44  &         5.43 \\
\textbf{WI($3,3;1$)} &         11.50  &         11.80  &         12.10  &         11.00  &         9.04  &         7.35 \\
\textbf{WI($3,5;1$)} & \textbf{12.00} &         13.00  & \textbf{12.80} &         11.60  &         9.18  & \textbf{7.96} \\
\textbf{WI($5,1;1$)} & \textbf{10.50} &         10.00  &          9.20  &          8.15  &         6.52  &         5.43 \\
\textbf{WI($5,3;1$)} &         11.00  & \textbf{13.20} &         12.40  &         10.90  &         8.72  &         7.09 \\
\textbf{WI($5,5;1$)} & \textbf{12.00} &         12.20  &         11.90  & \textbf{10.95} & \textbf{9.52} &         7.48 \\
\bottomrule
\end{tabular}
}
\end{minipage}
\hspace{0.02\textwidth}
\begin{minipage}[t]{0.47\textwidth}
\scalebox{0.75}{
\begin{tabular}{p{2.1cm}|rrrrrr}
\toprule
 \textbf{QN-WI} / $\Delta t$   & \hspace{0.27cm} 5.0 & \hspace{0.27cm} 2.0 & \hspace{0.27cm} 1.0 & \hspace{0.27cm} 0.5 & \hspace{0.27cm} 0.2 & \hspace{0.27cm} 0.1 \\
\midrule
\textbf{WI($1,1;1$)} & \textbf{10.50} & \textbf{10.00} &  \textbf{9.00} & \textbf{7.85} &         6.54  &         5.45 \\
\textbf{WI($1,3;1$)} &         11.50  &         10.60  &          9.70  &         8.85  &         7.42  &         6.60 \\
\textbf{WI($1,5;1$)} &         11.50  &         11.00  &          9.80  &         8.75  &         7.70  &         6.77 \\
\textbf{WI($3,1;1$)} & \textbf{10.50} & \textbf{10.00} &          9.20  &         8.10  & \textbf{6.44} & \textbf{5.43} \\
\textbf{WI($3,3;1$)} & \textbf{12.00} &         11.40  &         10.40  &         9.30  &         7.50  &         6.36 \\
\textbf{WI($3,5;1$)} & \textbf{12.00} & \textbf{11.80} &         11.10  &         9.85  & \textbf{8.16} & \textbf{6.89} \\
\textbf{WI($5,1;1$)} & \textbf{10.50} &         10.00  &          9.20  &         8.15  &         6.52  & \textbf{5.43} \\
\textbf{WI($5,3;1$)} &         11.50  & \textbf{11.80} & \textbf{11.30} & \textbf{9.85} &         8.00  &         6.82 \\
\textbf{WI($5,5;1$)} & \textbf{12.00} &         11.60  &         10.60  &         9.45  &         7.66  &         6.41 \\
\bottomrule
\end{tabular}
}
\end{minipage}
\vspace{1em}
\caption{Heat-heat coupling: average number of coupling iterations per time step for various coupling schemes and multi-rate setups \textbf{SC($\nD,\nN$)} and \textbf{WI($\nD,\nN;1$)}, respectively, and different window sizes $\Delta t$. As time integration method, implicit Euler is used for both solvers in all examples. The highest and the lowest number per column are always set in bold face.
}
\label{tab::qn-iterations}
\end{table}

\paragraph{Approximation Order in Time}

Above, we have examined the convergence speed of different iterative coupling algorithms and concluded that only quasi-Newton accelerated algorithms are able to reach a competitive number of iterations. In the following, we concentrate on the accuracy of the quasi-Newton based multi-rate setups \textbf{WI($\nD,\nN;p$)} and
\textbf{SC($\nD,\nN$)}. 

In a first step, we use $g_{\text{pol}}(t)$ and study for which setups and $p$ we can recover the exact solution. As stated in Sect.~\ref{sec:waveform}, a time stepping scheme of at least order $p$ and and waveform
interpolation of at least degree $p$ are required to achieve this. We consider the numerical solution $u$ of the partitioned heat equation to be equal to the exact solution $u_\text{exact}$, if the $L^2(\Omega)$ error is smaller than $10^{-12}$ in every time step. To avoid errors being introduced by the iterative coupling procedure, we use a very strict relative coupling tolerance of $10^{-12}$. As explained above, our discretization uses piecewise quadratic finite elements for the computation of $u$ to be able
to represent $u$ exactly in space.

Using \textbf{QN-WI} for \textbf{WI($\nD,\nN;p$)}, we were able to recover the exact solution for all multi-rate settings ($n_\mathcal{D}, n_\mathcal{N} \in \{1,2,3,5\}$) and all window sizes ($\Delta t \in \{0.0125, 0.025, 0.05, 0.1, 0.2, 0.5, 1.0\}$) for the following cases:
\begin{enumerate}
  \item $\gpolorder=1$, IE, and piecewise linear interpolation,
  \item $\gpolorder=1$, TR, and piecewise linear interpolation,
  \item \label{en:order2} $\gpolorder=2$, TR, and piecewise quadratic interpolation,
  \item \label{en:order3} $\gpolorder=3$, SDC, and piecwise cubic interpolation,
\end{enumerate}
For \ref{en:order2}, we also observed that, as expected, piecewise linear interpolation is not sufficient,
except for cases where interpolation is not required at all (such as \textbf{WI($5,5;1$)}). 
For \ref{en:order3}, we additionally tested that piecewise linear and piecewise quadratic interpolation
is never sufficient (not even for \textbf{WI($5,5;1$)}), since SDC internally performs two IE substeps, where a sufficiently accurate interpolation routine is required.

Using \textbf{QN-SC} for \textbf{SC($\nD,\nN$)}, we observed that only for the trivial non-multi-rate cases with $\nD=\nN=1$ using IE and $\gpolorder=1$ or TR and $\gpolorder=2$, where no interpolation is necessary, the exact solution can be recovered. The only multi-rate settings, where the exact solution could be obtained are the edge cases \textbf{SC($1,\nN$)}, $\nN \in \{1,2,3,5\}$ with $\gpolorder=1$ and IE. Here, the starting Dirichlet solver only requiring data at the end of the time window computes the correct flux for the end of the Neumann window, which is the only place where we measure the error. We expect the single time steps on Neumann side to not give the exact solution. For other multi-rate setups (such as \textbf{SC($5,2$)}), we could not obtain the exact solution, which was to be expected.

\begin{figure}[h!]
\begin{center}
\begin{tikzpicture}
\begin{loglogaxis}[
width=.6\textwidth,
height=.4\textwidth,
tuftelike,
legend cell align={left},
xlabel={$\Delta t$}, 
ylabel={$L^2(\Omega)$ error}, 
align=center, 
legend style={at={(1.1,0.5)},anchor=west, draw=none, align=left},
grid
]

\addplot[mark=x,color=red] table[x=dT,y=WR53, col sep=comma]{data/convergenceStudy/data_order/errors_raw_ie_l};
\addlegendentry{IE, \textbf{WI($5,3;1$)}};

\addplot[color=red, dashed] table[x=dT,y=WR53, col sep=comma]{data/convergenceStudy/data_order/errors_raw_ie_ss};
\addlegendentry{IE, \textbf{SC($5,3$)}};

\addplot[mark=square,color=blue] table[x=dT,y=WR53, col sep=comma]{data/convergenceStudy/data_order/errors_raw_tr_q};
\addlegendentry{TR, \textbf{WI($5,3;2$)}}

\addplot[color=blue, dashed] table[x=dT,y=WR53, col sep=comma]{data/convergenceStudy/data_order/errors_raw_tr_ss};
\addlegendentry{TR, \textbf{SC($5,3$)}};

\addplot[mark=*,color=black] table[x=dT,y=WR53, col sep=comma]{data/convergenceStudy/data_order/errors_raw_sdck16_c};
\addlegendentry{SDC, \textbf{WI($5,3;3$)}};

\draw[gray, thick] (axis cs:1,0.01) -- node[below]{$\mathcal{O}\left(\Delta t\right)$} (axis cs:0.01,0.0001);
\draw[gray, thick] (axis cs:1,0.0004) -- node[below=.2em]{$\mathcal{O}\left(\Delta t^2\right)$} (axis cs:0.01,0.00000004);
\draw[gray, thick] (axis cs:1,0.00004) -- node[below=.8em]{$\mathcal{O}\left(\Delta t^{3.5}\right)$} (axis cs:0.01,0.000000000004);

\end{loglogaxis}
\end{tikzpicture} \caption{Heat-heat coupling: convergence of the $L^2(\Omega)$-error at $T=1$ for the multi-rate setups \textbf{SC($5,3$)} and \textbf{WI($5,3;p$)} and various coupling, time intergration and interpolation methods.}
\label{fig:compare_MR53}
\end{center}
\end{figure}
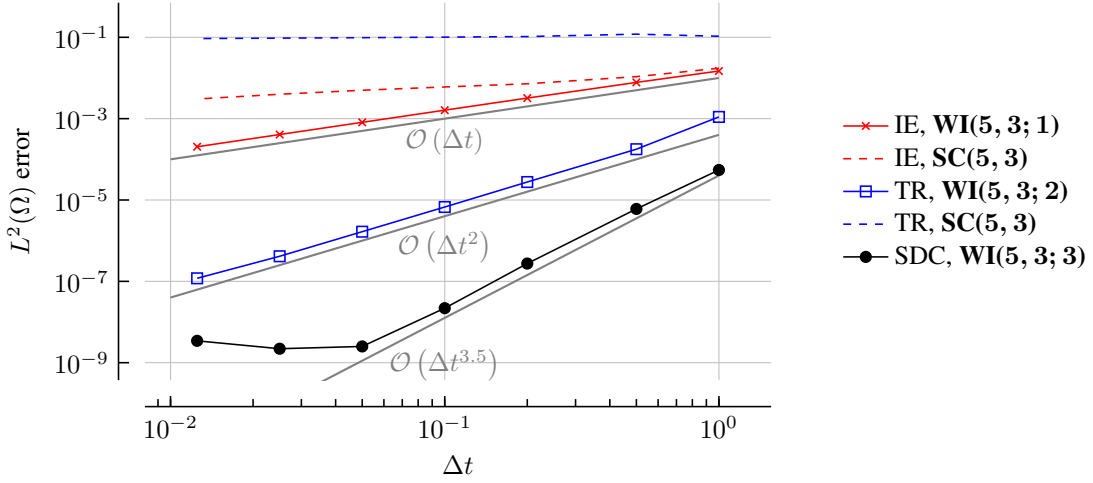

In a second step, we determine the convergence order of time stepping for the setup using $g_\text{tri}$.
For the convergence study, we determine the $L^2\left(\Omega\right)$ error of the numerical solution $u$ with respect to the exact solution $u_\text{exact}$. We consider the error at final simulation time $T=1$ and decrease the window size $\Delta t$. We use a moderate coupling tolerance of $10^{-5}$.
\autoref{fig:compare_MR53} shows a convergence study for \textbf{WI($5,3;p$)} and \textbf{SC($5,3$)}.
Note, that we intentionally
pick a setup where interpolated values are required by both participants of the coupled simulation at every time step. If \textbf{WI($5,3;p$)} is used for coupling in combination with a time stepping method of order $\ge p$, we expect convergence order $p$ as detailed in Sect.~\ref{sec:waveform} . We observe
this for IE and TR, combined with linear and quadratic interpolation, respectively. For the fourth-order SDC and cubic interpolation, we observe a convergence order of approximately $3.5$. Additionally, we observe that the accuracy of SDC is bounded by the coupling tolerance for small $\Delta t$ and the error does not scale with the time step anymore.
\textbf{SC($5,3$)}, on the other hand, gives slow to no convergence and significantly higher errors for IE and TR time integration methods.

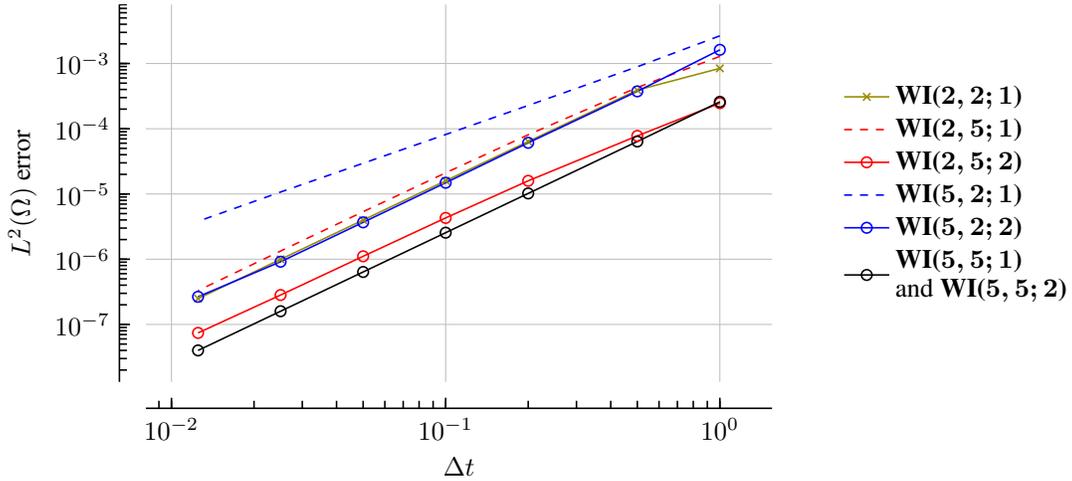
\begin{figure}[h!]
\begin{center}
\begin{tikzpicture}
\begin{loglogaxis}[
width=.6\textwidth,
height=.4\textwidth,
tuftelike,
legend cell align={left},
xlabel={$\Delta t$}, 
ylabel={$L^2(\Omega)$ error}, 
align=center, 
legend style={at={(1.1,0.5)},anchor=west, draw=none, align=left},
grid
]

\addplot[mark=x,color=olive] table[x=dT,y=WR22, col sep=comma]{data/convergenceStudy/data_tr_interp/errors_raw_tr_s_l};
\addlegendentry{\textbf{WI($2,2;1$)}};

\addplot[color=red, dashed] table[x=dT,y=WR25, col sep=comma]{data/convergenceStudy/data_tr_interp/errors_raw_tr_s_l};
\addlegendentry{\textbf{WI($2,5;1$)}};

\addplot[mark=o,color=red] table[x=dT,y=WR25, col sep=comma]{data/convergenceStudy/data_tr_interp/errors_raw_tr_s_q};
\addlegendentry{\textbf{WI($2,5;2$)}};

\addplot[color=blue, dashed] table[x=dT,y=WR52, col sep=comma]{data/convergenceStudy/data_tr_interp/errors_raw_tr_s_l};
\addlegendentry{\textbf{WI($5,2;1$)}};

\addplot[mark=o,color=blue] table[x=dT,y=WR52, col sep=comma]{data/convergenceStudy/data_tr_interp/errors_raw_tr_s_q};
\addlegendentry{\textbf{WI($5,2;2$)}};

\addplot[mark=o,color=black] table[x=dT,y=WR55, col sep=comma]{data/convergenceStudy/data_tr_interp/errors_raw_tr_s_q};
\addlegendentry{\textbf{WI($5,5;1$)}\\ and \textbf{WI($5,5;2$)}};

\end{loglogaxis}
\end{tikzpicture} \caption{Heat-heat coupling: Convergence of the $L^2(\Omega)$-error at $T=1$ for the trapezoidal rule as time integration method on both sides, \textbf{QN-WI}, and various multi-rate setups and interpolation methods.}
\label{fig:plot_tr_interp}
\end{center}
\end{figure}

To investigate the effect of multi-rate schemes being used, we perform another convergence study, where TR with \textbf{WI($2,2;p$)}, \textbf{WI($5,2;p$)}, \textbf{WI($2,5;p$)}, and \textbf{WI($5,5;p$)} is used (see \autoref{fig:plot_tr_interp}). The degree of the interpolating polynomials is irrelevant for cases, where no interpolation is needed, i.e., \textbf{WI($2,2;p$)} and \textbf{WI($5,5;p$)} reach second order regardless of the interpolation scheme being used. Note that this effect is only expected to occur for time stepping schemes such as TR or IE, where no additional function evaluations in time are required. For SDC, we did not observe this effect. For \textbf{WI($5,2;p$)} and \textbf{WI($2,5;p$)}, piecewise quadratic interpolation ($p = 2$) allows us to reach second order of the time stepping scheme, as expected. Here, the benefit of multi-rate is clearly visible: \textbf{WI($2,5;2$)} shows a significantly lower error than \textbf{WI($5,2;2$)} and almost the same as \textbf{WI($5,5;2$)}, which means that higher resolution in time for $\mathcal{N}$ is more beneficial than for $\mathcal{D}$. Presumably due to the term $x^2 \sin(t)$ in $u_{\text{exact}}$, which is significantly larger in $\Omega_\mathcal{N}$. 

If only piecewise linear interpolation($p = 1$) is applied, second order can be reached for certain cases (e.g. \textbf{WI($2,5;1$)}). This might be due to the smaller time steps on the Neumann side, from which the Dirichlet solver reads the interpolated data. In other words, it can be assumed that what we observe is the convergence of the dominating spatial error of the Dirichlet solver (using a $2.5$ times larger time step). \textbf{WI($5,2;1$)} shows the expected deterioration of the order below second order. We observe that, even though reaching second order convergence, the accuracy of \textbf{WI($2,5;1$)} is significantly lower than for the equivalent case with higher-order interpolation \textbf{WI($2,5;1$)}.

To summarize, we conclude that the converged solutions of \textbf{QN-WI} for \textbf{WI($\nD,\nN;p$)} fulfil the expectations in terms of accuracy. \textbf{SC($\nD,\nN$)}, on the other hand, yields only first order accuracy. We observe that interpolation plays an important role for waveform iteration. If we want to avoid order degradation, we have to use an interpolation with a sufficiently high polynomial degree. The required degree depends on the order of the time integration method, usually we require the polynomial degree $p$ to be at least equal to the order of the time stepping scheme. However, we notice that, for exceptional cases, $p$ lower than the order of the time stepping scheme can be sufficient.

\subsection{Fluid-Structure Interaction}
\label{ssec:FSI}

We now study whether we can generalize our observations from the relatively simple heat transfer problem to a more involved fluid-structure interaction problem. We show that \textbf{QN-WI} coupling allows us to couple two different second-order time integration methods, a trapezoidal rule used in the fluid solver and a Newmark-$\beta$ method used in the solid solver to an overall second-order coupled simulation even when combining different time step sizes in a multi-rate setup.   

\begin{figure}
\centering
\begin{minipage}{0.48\textwidth}
\includegraphics[width=0.95\textwidth]{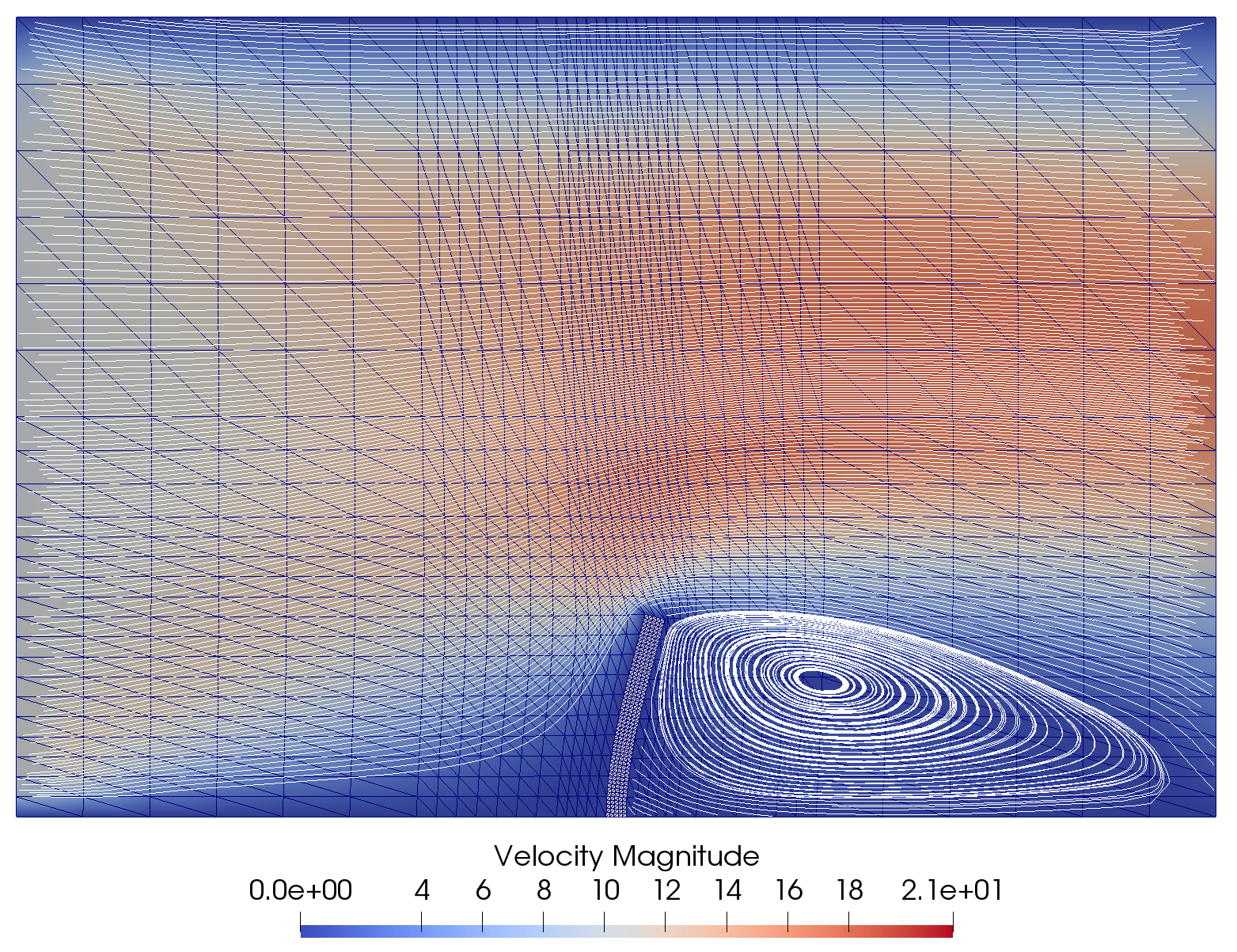}
\end{minipage}
\hspace{0.02\textwidth}
\begin{minipage}{0.48\textwidth}
\begin{tikzpicture}
\begin{axis}[
tuftelike,
width=.95 \textwidth,
height=.8 \textwidth,
xlabel={Time [s]}, 
ylabel={x-Displacement [m]}, 
align=center, 
    y tick label style={
        /pgf/number format/.cd,
            fixed,
            fixed zerofill,
            precision=2,
        /tikz/.cd
    },
grid,
legend style={at={(1,1)},anchor=north east},
ymax=0.25,
]

\addplot[color=black, dashed, thick] table[x=Time,y=Displacements10]{data/fsi_physics/precice-SOLIDZ-watchpoint-Tip-new.log};
\addlegendentry{Nutils-FEniCS};

\addplot[color=gray, thick] table[x=Time,y=Displacements00]{data/fsi_physics/precice-fenics-watchpoint-Tip.log};
\addlegendentry{OpenFOAM-FEniCS};

\end{axis}
\end{tikzpicture} %
\end{minipage}
\caption{
Perpendicular Flap. Left: snapshot of velocity magnitude and streamlines in the fluid domain and meshes of both domains at $t = 0.7\text{s}$ for a Nutils-FEniCS simulation. Right: tip displacement of the flap in $x$-direction over time, comparison between a OpenFOAM-FEniCS and a Nutils-FEniCS simulation for \textbf{QN-SC}, MR11, and $\Delta t=0.01$.
}
\label{fig::fsi}
\end{figure}
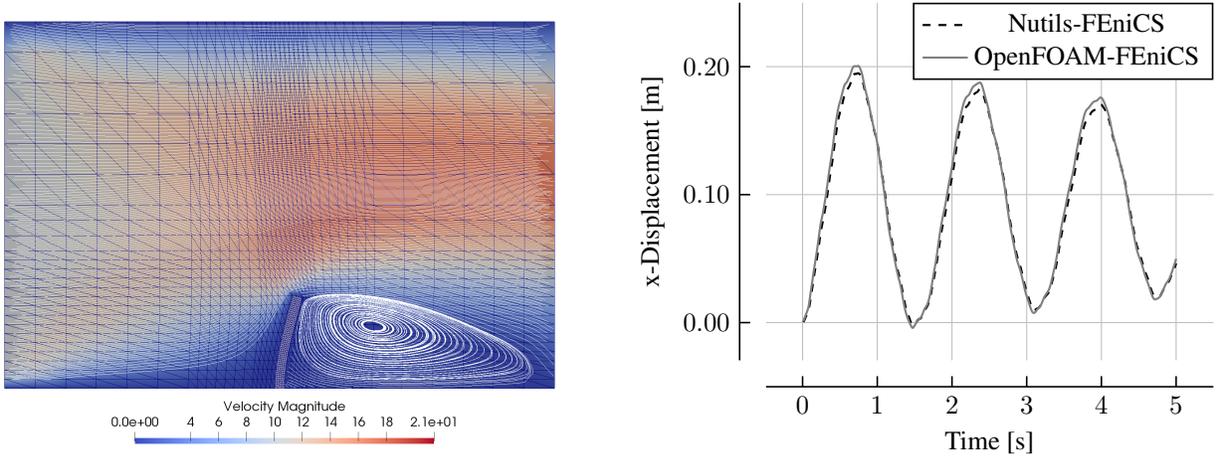

As FSI scenario, we study the Perpendicular Flap testcase from the preCICE tutorials\footnote{https://github.com/precice/tutorials}. An elastic flap with a width of $0.1\,\text{m}$ and a height of $1.0\,\text{m}$ is mounted in cross-flow at the center of a channel with a length of $6.0\,\text{m}$ and a height of $4.0\,\text{m}$. Figure \ref{fig::fsi} left shows the setup at maximal deflection.
This testcase is not experimentally validated, but cross-validated with different solver combinations. 
Compared to the preCICE tutorial, we slightly adjust the physical parameters: we use a higher fluid viscosity and a stiffer solid material to avoid requiring a finite element stabilization in the fluid solver and to use an only linear model in the structure solver, respectively. We set the fluid density to $\rho_F = 1.0 \,\text{kg}/\text{m}^3$, the kinematic viscosity to $\nu_F = 1.0 \,\text{m}^2/\text{s}$, the solid density to $\rho_S = 3.0 \cdot 10^3 \,\text{kg}/\text{m}^3$, the Young's modulus to $E = 4.0 \cdot 10^7 \,\text{kg}/\text{ms}^2$ and the Poisson ratio to $\nu=0.3$. On the left boundary a constant inflow profile in $x$-direction of $10\,\text{m/s}$ is prescribed. The right boundary is an outflow and the top and bottom of the channel as well as the surface of the flap are no-slip walls.     

For the solid domain, we use the linear elasticity solver developed in FEniCS and validated against CalculiX as part of the Bachelor thesis of Richard Hertrich\cite{Hertrich2019}, to which we also refer for details. We use a structured mesh with $5 \times 50$ cells, where each cell is cut into two triangles and $P_2$ finite elements for discretization. For time-integration, as mentioned previously, we use a second-order Newmark-$\beta$ method. The resulting linear system is solved by a direct solver.  

For the fluid domain, we implemented a new incompressible Navier-Stokes solver in Nutils specifically for this testcase. We use a structured adaptive mesh as depicted in Figure \ref{fig::fsi} left. We formulate the fluid equations fully-coupled and fully non-linearly in an arbitrary-Lagrangian-Eulerian framework and discretize with $Q^2-Q^1$ Taylor-Hood elements. The non-linearity is resolved by Newton's method with line-search up to an absolute residual of $10^{-9}$. A direct solver is used for the resulting linear systems. The mesh displacement is modeled by a simple Laplace equation, discretized by $Q^1$ elements and again solved with a direct solver. 
For time integration, we make use of the second-order trapezoidal rule, meaning we replace the time derivative of the velocity by the finite difference $(u_{i+1}-u_{i}) / \delta t$ and use a linear average between the variables of the current and the previous time step for the diffusive and the convective term. To avoid pressure oscillations, however, we take the pressure gradient and the continuity term only at the new time level.

To get to second order in an FSI setup, a few enhancements of the fluid time integration are necessary. Firstly, a second-order backward difference scheme is needed to compute the mesh velocity from the mesh displacement,
\[
u_{\text{mesh}} = \frac{3 \cdot d_{i+1} - 4 \cdot d_i + \cdot d_{i-1}}{2 \cdot \delta t} \,,
\]
where $u_{\text{mesh}}$ is the mesh velocity, $d$ the mesh displacement, and the index $i$ the time step. Secondly, the computation of forces acting on the solid needs to be constructed carefully. We compute fluxes over the fluid boundaries by evaluating the residual of the weak formulation at the current solution and extracting the coupling boundary\cite{Brummelen2012_FluxEvaluation}. Afterwards, we compute forces from the fluxes. To get second-order-accurate forces, we need to adapt the temporal discretization of two terms of the weak form compared to the trapezoidal form explained above. As time derivative of the velocity, we here use a second-order stencil, 
\[
\frac{3 \cdot u_{i+1}-4 \cdot u_{i} + u_{i-1}}{2 \cdot \delta t}\,,
\]     
where $u$ is the velocity. Furthermore, the trapezoidal rule with a fully-implicit pressure gradient as explained above only results in second-order pressure values at the midpoint $p_{i+0.5}$. Using these values at the time level $i+1$, as often done, is, however, only first-order accurate. To get second-order-accurate values at this time level, we extrapolate
\[
p_{i+1} = \frac{3}{2} \cdot p_{i+0.5} - \frac{1}{2} \cdot p_{i-0.5}
\]
in the pressure gradient term.

As the fluid solver and the solid solver use non-matching meshes at the coupling boundary, we require data mapping methods between both meshes. We, therefore, use radial-basis function interpolation with thin-plates splines as provided by preCICE\cite{Lindner2017}. In the fluid solver (Nutils), we directly read the displacement values at the Gauss points. FEniCS, however, requires a functional description of boundary conditions. We, therefore, read force values at the mesh vertices and construct the boundary function again by a radial-basis function interpolation; this time in the FEniCS adapter\cite{Hertrich2019}.    

To validate the new fluid solver, its mesh movement, and its force computation, we compare the Nutils-FEniCS results against OpenFOAM-FEniCS results from the preCICE tutorials. For both simulations, we use $\Delta t=0.01\,\text{s}$ and \textbf{QN-SC(1,1)} as coupling algorithm. Figure \ref{fig::fsi} shows a close agreement for the $x$-displacement of the tip until $T=5.0\,\text{s}$ for an initial fluid at rest.

\paragraph{Convergence Study}

To study the temporal convergence order of a coupled simulation, we compute the Perpendicular Flap testcase until $T = 0.01\,\text{s}$ with decreasing window sizes $\Delta t_l = 0.0025 \cdot 2^{-l}\,\text{s}$ for time levels $l=0,1,\ldots,4$. For all levels, we compute \textbf{QN-WI} with \textbf{WI($1,1;1$)}, \textbf{WI($2,3;2$)}, and \textbf{WI($3,2;2$)}. As reference solution for all methods, we use \textbf{WI($1,1;1$)} with time level $l=6$.
A crucial ingredient to achieve second order are suitable initial conditions. We start the fluid solver from an initial Stokes solution and then scale the forces on the structure until $T/2$ with a $C^1$-continuous $\sin^2(\pi (t+\delta t) / T)$ ramp. We measure the error in the displacement values at the coupling boundary to the reference solution $d_{\text{ref}}$. To smooth out variations, we consider the error over the complete time interval $[0,T]$,
\[
\vertiii{d_l-d_{\text{ref}}} := \left( \frac{1}{N_l}
\sum_{j=1}^{N_l} \left\| d_l(\Delta t_l \cdot j) - d_{\text{ref}}(\Delta t_l \cdot j) \right\|_{l^2(\Gamma)}^2
\right)^{1/2} \;\;\; l=0,1,\ldots, 4,
\]
where $N_l = \frac{T}{\Delta t_l}$ denotes the number of windows. For every time slice, we use 
\[
\| d \|_{l^2(\Gamma)} = \left( \sum_{\nu=1}^{m} d(x_{\nu})^2\right)^{1/2}\,,
\]  
where $x_{\nu}$ denotes the coordinates at the coupling interface and $m$ the number of vertices. Thus, $d(x_{\nu})$ are nodal evaluations.

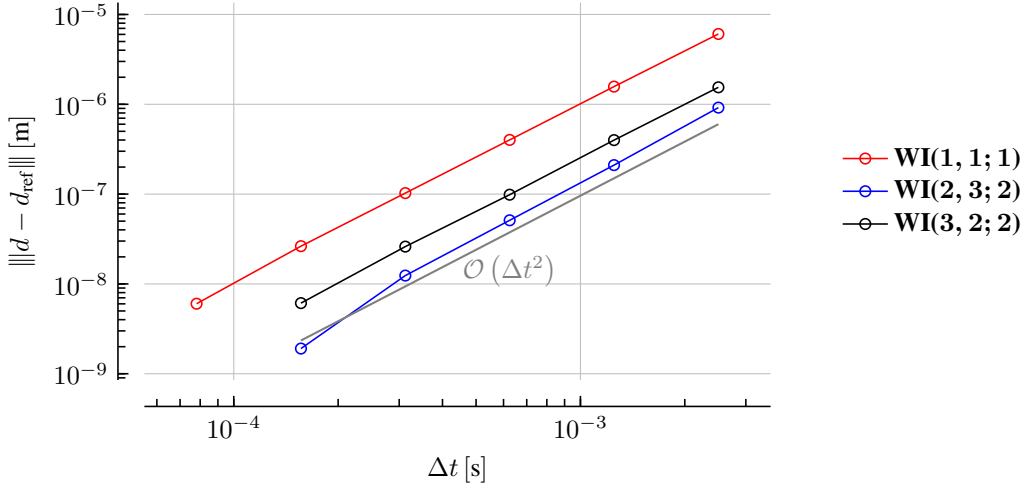
\begin{figure}[h!]
\begin{center}
\begin{tikzpicture}
\begin{loglogaxis}[
width=.6\textwidth,
height=.4\textwidth,
tuftelike,
legend cell align={left},
xlabel={$\Delta t\,[\text{s}]$}, 
ylabel={$\vertiii{d-d_{\text{ref}}}\,[\text{m}]$}, 
align=center, 
legend style={at={(1.1,0.5)},anchor=west, draw=none, align=left},
grid
]

\addplot[mark=o,color=red] table[x=dT,y=MR11, col sep=comma]{data/fsi/errors_raw};
\addlegendentry{\textbf{WI($1,1;1$)}};

\addplot[mark=o,color=blue] table[x=dT,y=MR23, col sep=comma]{data/fsi/errors_raw};
\addlegendentry{\textbf{WI($2,3;2$)}};

\addplot[mark=o,color=black] table[x=dT,y=MR32, col sep=comma]{data/fsi/errors_raw};
\addlegendentry{\textbf{WI($3,2;2$)}};

\draw[gray, thick] (axis cs:0.0025,0.0000006) -- node[below=.5em]{$\mathcal{O}\left(\Delta t^2\right)$} (axis cs:0.00015625,0.00000000234375);

\end{loglogaxis}
\end{tikzpicture}
 \caption{
\label{fig:fsi_order}
Perpendicular Flap: Convergence of the displacement values at the coupling interface towards the reference simulation over the complete time interval $[0,0.01]\,\text{s}$. 
Three multi-rate setups for \textbf{QN-WI} are compared.}
\end{center}
\end{figure}

Figure \ref{fig:fsi_order} compares the errors for all time levels. QNWI coupling achieves a clear second-order convergence for all three considered multi-rate setups. The errors of \textbf{WI($1,1;1$)} and \textbf{WI($3,2;2$)} appear to be shifted by one time level, which indicates that the solid time step size dominates the overall error. This again shows the great potential of multi-rate time stepping as the fluid problem is the far more expensive one. Thus, using a smaller time step size for the structure solver is almost for free. We compute the coupling iterations up to a relative residual error of $10^{-6}$ and list the average number of iterations in Tab.~\ref{tab:fsi}. We can observe a slight increase of required iterations with decreasing window size, contrary to our observations for the heat transfer testcase in the last section. This does not come as a surprise as, for incompressible FSI and vanishing structural mass, a similar behavior is reported in literature\cite{VanBrummelen2009_AddedMassEffect}. It opens, however, the possibility for optimization as increasing the window size while keeping the time step size constant might lead to a more efficient numerical method, besides also having all the computing advantages of increased asynchronicity (such as reduced communication and better balancing of varying computational load). A detailed study is, however, beyond the scope of this paper. Lastly, comparing the required iterations of \textbf{WI($2,3;2$)} with \textbf{WI($3,2;2$)} could indicate that the structural time step size is the critical part.

\begin{table}[h!]
\centering
\begin{tabular}{r|ccccc}
\toprule
$\Delta t_l [\text{s}]$ & $0.0025 \cdot 2^0$ & $0.0025 \cdot 2^{-1}$ & $0.0025 \cdot 2^{-2}$ & $0.0025 \cdot 2^{-3}$ & $0.0025 \cdot 2^{-4}$ \\
\midrule
\textbf{WI($1,1;1$)} & 4.00 & 4.50 & 4.81 & 5.50& 5.64 \\
\textbf{WI($2,3;2$)} & 5.25 & 5.63 & 6.57 & 7.31 & 7.42 \\
\textbf{WI($3,2;2$)} & 4.50& 4.75 & 5.31 & 5.63 & 6.33 \\
\bottomrule
\end{tabular}
\caption{\label{tab:fsi}
Perpendicular Flap: Average number of quasi-Newton iterations per time window until $T=0.01\,\text{s}$ for \textbf{QN-WI}, various multi-rate setups, and a decreasing window size from left to right.}
\end{table}
\section{Conclusions and Outlook}\label{sec::conclusions}

This contribution shows that using waveform iteration in a black-box approach for partitioned multi-physics simulations is possible. Black-box means here in particular that the required interpolation in time can be computed without using discretization details of the coupled codes. An elegant prototype software design using a middle-layer called \textit{waveform bindings} allows us to test arbitrary mutli-rate settings with the coupling library preCICE without any alterations to the latter.
Arbitrary multi-rate settings means that we allow any (and different) numbers of time steps in the involved solvers within a so-called coupling window. Our results show that (i) the implementation actually achieves the predicted order of consistency in time and that (ii) we can generalize interface quasi-Newton methods to the space-time domain.
For more complex scenarios such as fluid-structure interaction, we observed that actually achieving higher order requires careful consideration of all solver components, of mesh moving, and of evaluations of coupling variables.

Our current implementation has a couple of limitations, which we intend to tackle in future work: 
(i) We have to move the functionality of the waveform bindings to preCICE and decide how to do this while maintaining backwards compatibility of the application user interface as far as possible. This also removes the technical restriction to Python. (ii) Features that optimize the performance of quasi-Newton schemes need to be provided also for waveform iterations. This holds in particular for re-use of previous time step data and improved initial guesses for the initial approximation (currently constant extrapolation) of interface data before the first iteration in a new time window. (iii) Our methods should be generalized to adaptive time stepping, which is less of a mathematically unsolved problem, but a software issue as it might require varying numbers of time steps and, thus, varying numbers of coupling data structures in preCICE from one time window to the next. 
In addition to these enhancements, further studies for arbitrary high order (as in \cite{Huang2019}) will be necessary.

\section*{Acknowledgments}
We want to thank Harald van Brummelen for the fruitful discussions concerning the results of our research, and in particular concerning the second-order accuracy of the fluid solver in Sect.~\ref{ssec:FSI}, Gertjan van Zwieten for the high-quality Nutils support and Richard Hertrich for the implementation and validation of the structural solver in FEniCS.

Furthermore, we thankfully acknowledge the funding from the European Union's Horizon 2020 research and innovation program under the Marie Sklodowska-Curie grant agreement No 754462 as well as funding from SPPEXA, the German Science Foundation Priority Program 1648 -- Software for Exascale Computing.

\bibliographystyle{unsrt}
\bibliography{literature}%

\clearpage

\end{document}